\def\mapright#1{\smash{\mathop{\longrightarrow}\limits^{#1}}}

\def\mapup#1{\uparrow\rlap{$\vcenter{\hbox{$\scriptstyle#1$}}$}}
\def\mapupleft#1{\uparrow\llap
                   {$\vcenter{\hbox{$\scriptstyle#1$}}\;\;$}}

\documentclass{amsart}
\usepackage{amscd, amsmath, amsthm, amssymb}

\newtheorem{theorem}{Theorem}[section]
\newtheorem{lemma}[theorem]{Lemma}
\newtheorem{proposition}[theorem]{Proposition}

\theoremstyle{definition}     
\newtheorem{definition}[theorem]{Definition}

\newtheorem{claim}[theorem]{Claim}

\theoremstyle{remark}
\newtheorem{remark}[theorem]{Remark}

\numberwithin{equation}{section}

\begin{document}

{\bf Version: 26 November 2003}
\par \vskip 3pc

\title[Alternating group, Leech lattice and K3 surfaces]
{The alternating group of degree 6 in geometry of the Leech
lattice and K3 surfaces}

\author[J.H. Keum]{JongHae Keum}
\address{School of Mathematics, Korea Institute for Advanced Study, Dongdaemun-gu, Seoul 130-722, Korea}
\email{jhkeum@kias.re.kr}

\author[K. Oguiso]{Keiji Oguiso}
\address{Graduate School of Mathematical Sciences,
University of Tokyo, Komaba, Meguro-ku,
Tokyo 153-8914, Japan
}
\email{oguiso@ms.u-tokyo.ac.jp}

\author[D. -Q. Zhang]{De-Qi Zhang}
\address{Department of Mathematics, National University of Singapore,
2 Science Drive 2, Singapore 117543, Singapore}
\email{matzdq@math.nus.edu.sg}


\begin{abstract}
The alternating group of degree $6$ is located at the junction
of three series of simple non-commutative groups : simple sporadic groups, alternating groups and simple groups of Lie type. It plays a very special role in the theory of finite groups. We shall study its
new roles both in a finite geometry of certain pentagon in the Leech lattice and also in a complex algebraic geometry of $K3$ surfaces.
\end{abstract}

\maketitle


\setcounter{section}{0}
\section{ Introduction}
The alternating group $A_{6}$ plays a very special role in
the theory of finite groups [Su]. Indeed, $A_{6}$ is a simple
group which is located at a sort of junction
of three series of simple non-commutative groups : simple sporadic groups,
alternating groups and simple groups of Lie type. More explicitly, there are
distinguished isomorphisms $[M_{10}, M_{10}] \simeq
A_{6} \simeq  \text{PSL}(2, 9)$; see for instance Conway-Sloane [CS, Chapter 10].
Though $M_{10}$, the Mathieu group of degree $10$, is not itself
a simple group, it falls into a sequence  $M_{10} < M_{11} < M_{12}$
of maximal subgroups of the smallest sporadic simple groups $M_{11}$ and
$M_{12}$.

From a slightly different view, in contrast to the fact that
$\text{Aut}(A_{n}) \simeq S_{n}$, whence $\text{Out}(A_{n}) \simeq
C_{2}$ when $n \not = 6$ and $n \geq 3$, the outer automorphism
group $\text{Out}(A_{6})$ is isomorphic to a bigger group
$C_{2}^{\oplus 2}$. Corresponding to the three involutions,
$\text{Aut}(A_{6})$ has three index $2$ subgroups $A_{6} < G <
\text{Aut}(A_{6})$, which are $S_{6}$, $\text{PGL}(2, 9)$ and
$M_{10}$. According to Suzuki [Su, Page 300], it is this
extraordinary property which seems to make the classification of
simple groups deep and difficult. This property of $A_{6}$ also
plays a crucial role in our note (the proof of Proposition 2.6).

The aim of this note is to study roles played by $A_{6}$ both in a finite geometry in the Leech lattice $\Lambda$ (or in a slightly different language, in the set of Leech roots of $II_{1, 25} := \Lambda \oplus U$) and also in a complex algebraic geometry of $K3$ surfaces (Theorems 2.3, 3.1, 5.1 and Proposition 3.5).

We first show that $A_{6}$ can be characterized as the pointwise
stabilizer group of some uniquely determined pentagon in the Leech
lattice, or equivalently, a special configuration of Leech roots
of Coxeter-Dynkin type $A_{2}^{\oplus 2} \oplus A_{1}^{\oplus 2}$
(Theorem 2.3). This is an analogue of results by Curtis for
$S$-lattices [Cu] and by  Finkelstein for some maximal subgroups
of the Conway group $\cdot 3$ [Fi] and a table in
[CS, Page 291]. On the Leech lattice, Leech roots and the Conway
groups, we refer the readers to the standard reference book [CS,
Chapter 10, 28]; see also Section 2 below for a brief summary. We
note that our $A_{6}$ is in the subgroup $\cdot 3$ of $\cdot 0 (:=
O(\Lambda))$ but not in the Mathieu subgroups $M_{22} < M_{24}$
embedded in $\cdot 0$ in a standard way [Cu].

We then apply this characterization in our study of group symmetries on $K3$ surfaces. By definition,
a $K3$ surface is a simply-connected compact complex surface admitting a nowhere vanishing global holomorphic $2$-form.
They form $20$-dimensional moduli;
for more details about $K3$ surfaces, see for instance [BPV]. According to Mukai [Mu, Main theorem],
there are exactly 11 maximal finite groups each of which acts on some $K3$ surface symplectically. They all can be embedded
into the Mathieu group $M_{23}$. Among 11 such groups,
simple groups are only $\text{PSL}(2, 7)$ and the present
$A_{6} = \text{PSL}(2, 9)$.

Our goal is to show the existence and uniqueness of the triplet
$(F, \tilde{A}_{6}, \rho_{F})$ of  a $K3$ surface $F$ and its
finite group action $\rho_{F} : \tilde{A}_{6} \times F
\longrightarrow F$ of $\tilde{A}_{6}$ on $F$, up to isomorphisms
(Theorems 3.1 and 5.1). Here the group $\tilde{A}_{6}$ is an
extension of $A_{6}$ by $\mu_{4} \simeq C_{4}$, which turns out to
be the unique maximal possible finite extension of $A_{6}$ in the
automorphism groups of $K3$ surfaces (Theorem 5.1 and Proposition 4.1). 
This part is much
inspired
by the work of Kondo [Ko2, Ko3]. We remark that a work for the other
simple group $\text{PSL}(2, 7)$ has been carried out in [OZ2].

Our $K3$ surface $F$ turns out to be isomorphic to the minimal
resolution of the branched double cover of the elliptic modular
surface with level $3$ structure. We will write down explicitly
the bidegree $(2, 3)$ equation for a canonical model of $F$ in
$\mathbf P^{1} \times \mathbf P^{2}$ (Proposition 3.5). On the
other hand, we also see that $\text{Pic}(F)^{\tilde{A}_{6}} =
\mathbf Z H_{F}$ and $(H_{F}^{2}) = 20$ (Propositions 4.1 and
4.5). So the action of $\tilde{A}_{6}$ on $F$ is not induced by
$\text{PGL}(\mathbf P^{1}) \times \text{PGL}(\mathbf P^{2})$. The
invariant degree $(H_{F}^{2}) = 20$ also tells that our example
$(F, A_{6})$ ($A_{6} < \tilde{A}_{6}$) is not isomorphic to
Mukai's example $(X_{6}, A_{6})$ of polarized $K3$ surface of
degree $6$ with symplectic group action of $A_{6}$ in [Mu, Example
0.4, No. 2]. One can check that in Mukai's example the maximal
extension of $A_6$ in the full automorphism group
$\text{Aut}(X_6)$ is the symmetric group $S_6$.

We also remark that one can construct a smooth non-isotrivial
family of $K3$ surfaces $f : \mathcal X \longrightarrow \mathbf
P^{1}$ such that the fibres $\mathcal X_{t}$ admit $A_{6}$-actions
in exactly the same manner as in [OZ2, Appendix].

It would be very interesting to see the full automorphism group $\text{Aut}(F)$, which is of infinite order [SI].
In this direction, readers may refer to
[Vi], [Ko1], [KK] and [DK] for other $K3$ surfaces.

\par \vskip 1pc

{\it Acknowledgement.} Part of this work was carried out during
the first two authors' stay at National University of Singapore in
September 2003. They would like to express their thanks to the
Department of Mathematics and staff members there for financial
support and warm hospitalities. The second author would also like
to express his thanks to Professors A. A. Ivanov and A. Matsuo for
valuable comments.

\section{Uniqueness of the Leech roots of Coxeter-Dynkin type
$A_{2}^{\oplus 2}\oplus A_{1}^{\oplus 2}$}

First we briefly recall some necessary notations and facts about Leech lattice from [CS, Chapter 10], [Ko1] and [DK]. Let

$$\Omega := \mathbf P^1(\mathbf F_{23}) = \{\infty\, ,\, 0\, ,\, 1\, ,\, \cdots \, ,\, 22\}\,\, .$$

We denote by $2^{\Omega}$ the power set of $\Omega$. $2^{\Omega}$
has a structure of $24$-dimensional vector space over $\mathbf
F_{2}$, in which the sum is defined to be the symmetric
difference. Let $N := \Omega - \{a^{2} \vert a \in \mathbf
F_{23}\}$, and $N_{\infty} := \Omega$, $N_{i} := \{n -i \vert n
\in N\}$ ($i \in \mathbf F_{23}$). The subspace $\mathcal C$
spanned by $N_{i}$ ($i \in \Omega)$ is called the binary Golay
code. We call an element of $\mathcal C$ a $\mathcal C$-set. Let
$K$ be a $\mathcal C$-set, i.e. $K \in \mathcal C$. Then as a
subset of $\Omega$, $\vert K \vert$ is either $0$, $8$, $12$, $16$
or $24$. An element $K \in \mathcal C$ is called an octad (resp. a
dodecad) if $\vert K \vert = 8$ (resp. if $\vert K \vert = 12$).
It is well-known that the set of octads forms the so-called
Steiner system $St(5, 8, 24)$ of $\Omega$. There are exactly $759$
octads and they are explicitly listed in [To].

Next consider the $24$-dimensional Euclidean space $\mathbf
R^{24}$ with orthonormal basis $\langle \mu_{i} \rangle_{i \in
\Omega}$, i.e. $\mu_{i}.\mu_{j} = \delta_{ij}$ for $i, j \in
\Omega$. For any subset $A$ of $\Omega$, let $\mu_{A}$ denote the
vector $\sum_{i \in A} \mu_{i}$. Let $\Lambda$ be the $\mathbf
Z$-submodule of $\mathbf R^{24}$ spanned by the vectors $2\mu_{K}$
and $\mu_{\Omega} - 4\mu_{\infty}$, where $K$ runs through all
octads. Define the bilinear form on $\Lambda$ by $(U, V) :=
-U.V/8$ for $U, V \in \Lambda$. Then it is well known that $(U, V)
\in \mathbf Z$ and $\Lambda := (\Lambda, (*, *))$ forms an even
unimodular negative definite lattice of rank $24$. To be precise,
this is the so called {\it Leech lattice}. It is also well known that
$\Lambda$ contains no element $V$ with $(V^{2}) = -2$.

The following Theorem (see for instance [CS, Chapter 10, Theorem 25]) gives us a more concrete picture of the Leech lattice:

\begin{theorem}\label{theorem:leechv} The vector
$\sum_{i \in \Omega}x_{i}\mu_{i}$ ($x_{i} \in \mathbf Z$) is in $\Lambda$
if and only if the following four conditions are satisfied:
\begin{list}{}{
\setlength{\leftmargin}{10pt}
\setlength{\labelwidth}{6pt}
}
\item[(i)]
$x_{\infty} \equiv x_{0} \equiv x_{1} \equiv \cdots \equiv x_{21}
\equiv x_{22}\, (\text{mod}\, 2)$; denote by $m \equiv x_{\infty}$
$(\text{mod} \, 2)$.
\item[(ii)]
For each $j \in \{0, 1, 2, 3\}$, the set $\{i \in \Omega\, \vert\, x_{i} \equiv j\, (\text{mod}\, 4)\}$ is a $\mathcal C$-set.
\item[(iii)] $\sum_{i \in \Omega} x_{i} \equiv 4m$ $(\text{mod}\, 8)$.
\end{list}
\end{theorem}
For instance, from this theorem, one knows that every
element $V$ with $(V^{2}) = -4$ has one of the forms: $((\pm 2)^{8},
0^{16})$, where the non-zero coordinates have positive product and
are in the place of an octad; $(\mp 3, (\pm 1)^{23})$, where the
lower signs are taken on a $\mathcal C$-set; and $((\pm 4)^{2},
0^{22})$ with no extra condition.

The orthogonal group of the Leech lattice is denoted by $\cdot 0$.
It is well known that $\cdot 0$ acts on the set of vectors
$(V^{2}) = -2m$ ($m = 2$ or $3$) transitively. We denote the
stabilizer group of a vector $V$ with $(V^{2}) = -2m$ by $\cdot
m$. For more details about the Leech lattice, see for instance
[CS, Chapter 10].

Let $\Lambda$ be the Leech lattice and $\Pi = II_{1, 25} := \Lambda
\oplus U$ the unique even unimodular lattice of index $(1,
25)$. Here $U$ is the hyperbolic plane, i.e. the lattice $\mathbf
Z^{2}$ equipped with a bilinear form $((l_{1}, m_{1}),(l_{2},
m_{2}) ) = l_{1}m_{2} + l_{2}m_{1}$. This $U$ is the unique even
unimodular hyperbolic lattice. Let $w := (0, 0, 1) \in \Pi$ be the
Weyl vector. Set
$$\Pi_{2} := \{r \in \Pi \,\vert\, (r^{2}) = -2, \,\,
(r, w) = 1\}.$$
An element of $\Pi_{2}$ is called a Leech root of
$\Pi$. The positive cone, denoted by $\mathcal P$, is the one of
the two connected components of $\{v \in \Pi \otimes \mathbf R
\vert (v^{2}) > 0\}$ whose closure contains $w$. We set
$$\mathcal D = \{v \in \mathcal P\, \vert \, (v, r) \geq 0, \quad
\forall r \in \Pi_{2} \}.$$
It is well known that the correspondence
$$\Lambda \ni X \leftrightarrow x := (X, 1, -\frac{X^{2}}{2} -1)
\in \Pi_{2}$$ between $\Lambda$ and $\Pi_{2}$
is bijective, and under this bijection, we have a natural identification $\cdot \infty  = \text{Aut}(\mathcal D)$ [CS, Chapter 27, Theorem 1].
Here $\cdot \infty = \Lambda : \cdot 0$ is the group of affine isometries (namely including translations) of $\Lambda$. We also note that the Weyl vector $w$ is stable under $\text{Aut}(\mathcal D)$, i.e. $\varphi(w) = w$ for each $\varphi \in \text{Aut}(\mathcal D)$.

Let us consider the following $6$ vectors in $\Lambda$:
$$C = 4\nu_{\infty} + \nu_{\Omega}\,,\quad
Z = 0\,,\quad X_0= 4\nu_{\infty} + 4\nu_{0}\,,$$
$$R_{0} = 2\nu_{K_{0}}\,,\quad X_{1} = 2\nu_{K_{1}}\,,\quad X_{2} = 2\nu_{K_{2}}\,,$$
where
$$K_0=\{\infty, 1,2,3,4,6,15,18\}\,,$$
$$K_1=\{\infty, 0,1,2,3,5,14,17\}\,,\quad K_2=\{\infty, 0,1,2,4,13,16,22\}\,.$$
Let
$$c = (C, 1, 2)\,, \quad z = (0, 1, -1)\,,\quad x_{0} = (X_{0}, 1, 1)\,,$$
 $$r_{0} = (R_0, 1,1)\,, \quad x_{1} = (X_{1}, 1, 1)\,,\quad x_{2} = (X_{2}, 1, 1)$$
be the corresponding Leech roots. We employ here the same notation
as in [DK], but we rename $X$ there as $C$ here. Set:

$$\mathcal R := \{c, z, x_{0}, r_{0}, x_{1}, x_{2}\}\,\, .$$

It is easy to verify that $\mathcal R$ forms a Coxeter-Dynkin diagram of type $A_{2}^{\oplus 2}\oplus A_{1}^{\oplus 2}$:

\par \vskip 3pc

\setlength{\unitlength}{3947sp}%
\begingroup\makeatletter\ifx\SetFigFont\undefined%
\gdef\SetFigFont#1#2#3#4#5{%
  \reset@font\fontsize{#1}{#2pt}%
  \fontfamily{#3}\fontseries{#4}\fontshape{#5}%
  \selectfont}%
\fi\endgroup%
\begin{picture}(5633,1341)(826,-1844)
{ \thinlines
\put(2176,-586){\circle{150}}
}%
{ \put(3301,-586){\circle{150}}
}%
{ \put(4351,-586){\circle{150}}
}%
{ \put(5401,-586){\circle{150}}
}%
{ \put(6376,-586){\circle{150}}
}%
{ \put(1051,-586){\circle{150}}
}%
{ \put(3376,-586){\line( 1, 0){900}}
}%
{ \put(1201,-586){\line( 1, 0){900}}
}%
\put(826,-1036){\makebox(0,0)[lb]{\smash{\SetFigFont{17}{20.4}{\rmdefault}{\mddefault}{\itdefault}{ c}%
}}}
\put(1876,-1036){\makebox(0,0)[lb]{\smash{\SetFigFont{17}{20.4}{\rmdefault}{\mddefault}{\itdefault}{ z}%
}}}
\put(3076,-1036){\makebox(0,0)[lb]{\smash{\SetFigFont{17}{20.4}{\rmdefault}{\mddefault}{\itdefault}{ r}%
}}}
\put(3226,-1111){\makebox(0,0)[lb]{\smash{\SetFigFont{14}{16.8}{\rmdefault}{\mddefault}{\updefault}{ 0}%
}}}
\put(4126,-1036){\makebox(0,0)[lb]{\smash{\SetFigFont{17}{20.4}{\rmdefault}{\mddefault}{\itdefault}{ x}%
}}}
\put(4351,-1111){\makebox(0,0)[lb]{\smash{\SetFigFont{14}{16.8}{\rmdefault}{\mddefault}{\updefault}{ 0}%
}}}
\put(5176,-1036){\makebox(0,0)[lb]{\smash{\SetFigFont{17}{20.4}{\rmdefault}{\mddefault}{\itdefault}{ x}%
}}}
\put(5401,-1111){\makebox(0,0)[lb]{\smash{\SetFigFont{14}{16.8}{\rmdefault}{\mddefault}{\updefault}{ 1}%
}}}
\put(6151,-1036){\makebox(0,0)[lb]{\smash{\SetFigFont{17}{20.4}{\rmdefault}{\mddefault}{\itdefault}{ x}%
}}}
\put(6376,-1111){\makebox(0,0)[lb]{\smash{\SetFigFont{14}{16.8}{\rmdefault}{\mddefault}{\updefault}{ 2}%
}}}
\put(3376,-1786){\makebox(0,0)[lb]{\smash{\SetFigFont{14}{16.8}{\rmdefault}{\mddefault}{\updefault}{ Figure 1}%
}}}
\end{picture}
\par \vskip 1pc
We first remark the following easy but important fact. We recall a
few necessary notations. Let $(M, \langle*, *\rangle)$ be an even
non-degenerate lattice and $M^{*} = \text{Hom}\,(M, \mathbf Z)$
its dual. We embed $M \subseteq M^{*}$ by the non-degeneracy of
$\langle*,*\rangle$, and define $A_{M} := M^{*}/M$. Since $M$ is
even, we have a natural quadratic form $q_{M} : A_{M}
\longrightarrow \mathbf Q/ 2\mathbf Z$ given by: $v\, \text{mod}\,
M \mapsto v^{2}\, \text{mod}\, 2\mathbf Z$.

\begin{lemma}\label{lemma:R} Let $\mathcal R' = \{c', z', x_{0}', r_{0}', x_{1}', x_{2}'\} \subset \Pi_{2}$ be a set consisting of six Leech roots
forming a Coxeter-Dynkin diagram of type $A_{2}^{\oplus 2}\oplus
A_{1}^{\oplus 2}$. Let $R'$ be the sublattice of $\Pi$ generated
by $\mathcal R'$. Then $R'$ is isomorphic to $A_{2}^{\oplus
2}\oplus A_{1}^{\oplus 2}$ and $R'$ is a primitive sublattice of
$\Pi$.
\end{lemma}
\begin{proof}
Note that the discriminant group $A_{R'}\simeq (\mathbf Z/3)^2
\oplus (\mathbf Z/2)^2$ is generated by $(c'+2z')/3$,
$(r_{0}'+2x_{0}')/3$, $x_{1}'/2$, $x_{2}'/2$. A direct calculation
shows that $A_{R'}$ contains no isotropic element, whence
$R'^*\cap\Pi =R'$. This completes the proof.
\end{proof}
Note also that the five vectors $C, X_0, R_0, X_1, X_2\in
\Lambda$, which are projections of the Leech roots $c, x_0, r_0,
x_1, x_2$, form the following pentagon, where the number on an
edge is the intersection number of vectors joined by the edge,
e.g. $(R_0,\, X_0)=-1$, $C^2=-6$ (and the meaning of numbers is different from
that in [Cu]):

\par \vskip 3pc

\setlength{\unitlength}{3947sp}%
\begingroup\makeatletter\ifx\SetFigFont\undefined%
\gdef\SetFigFont#1#2#3#4#5{%
  \reset@font\fontsize{#1}{#2pt}%
  \fontfamily{#3}\fontseries{#4}\fontshape{#5}%
  \selectfont}%
\fi\endgroup%
\begin{picture}(5700,5713)(3451,-5144)
{ \thinlines
\put(8926,-1411){\circle{150}}
}%
{ \put(5101,-3811){\circle{150}}
}%
{ \put(6376,164){\circle{150}}
}%
{ \put(3976,-1336){\circle{150}}
}%
{ \put(7951,-3811){\circle{150}}
}%
{ \put(6301,164){\line(-3,-2){2215.385}}
}%
{ \put(6451,164){\line( 5,-3){2426.471}}
}%
{ \put(3976,-1411){\line( 1,-2){1140}}
}%
{ \put(5176,-3811){\line( 1, 0){2700}}
}%
{ \put(8926,-1486){\line(-2,-5){900}}
}%
{ \put(6451, 89){\line( 2,-5){1525.862}}
}%
{ \put(4051,-1336){\line( 1, 0){4800}}
}%
{ \put(4051,-1411){\line( 5,-3){3838.235}}
}%
{ \put(5176,-3736){\line( 5, 3){3694.853}}
}%
{ \put(5176,-3661){\line( 1, 3){1237.500}}
}%
\put(4126,-2686){\makebox(0,0)[lb]{\smash{\SetFigFont{17}{20.4}{\rmdefault}{\bfdefault}{\updefault}{ -2}%
}}}
\put(7651,-586){\makebox(0,0)[lb]{\smash{\SetFigFont{17}{20.4}{\rmdefault}{\bfdefault}{\updefault}{ -3}%
}}}
\put(4726,-586){\makebox(0,0)[lb]{\smash{\SetFigFont{17}{20.4}{\rmdefault}{\bfdefault}{\updefault}{ -3}%
}}}
\put(5851,-2536){\makebox(0,0)[lb]{\smash{\SetFigFont{17}{20.4}{\rmdefault}{\bfdefault}{\updefault}{ -2}%
}}}
\put(6226,-4036){\makebox(0,0)[lb]{\smash{\SetFigFont{17}{20.4}{\rmdefault}{\bfdefault}{\updefault}{ -2}%
}}}
\put(5776,-1861){\makebox(0,0)[lb]{\smash{\SetFigFont{17}{20.4}{\rmdefault}{\bfdefault}{\updefault}{ -3}%
}}}
\put(3451,-1111){\makebox(0,0)[lb]{\smash{\SetFigFont{17}{20.4}{\rmdefault}{\mddefault}{\updefault}{     (-4)}%
}}}
\put(8476,-2761){\makebox(0,0)[lb]{\smash{\SetFigFont{17}{20.4}{\rmdefault}{\bfdefault}{\updefault}{ -2}%
}}}
\put(4651,-4186){\makebox(0,0)[lb]{\smash{\SetFigFont{17}{20.4}{\rmdefault}{\mddefault}{\updefault}{ X}%
}}}
\put(4876,-4336){\makebox(0,0)[lb]{\smash{\SetFigFont{14}{16.8}{\rmdefault}{\bfdefault}{\updefault}{ 1}%
}}}
\put(3676,-1636){\makebox(0,0)[lb]{\smash{\SetFigFont{14}{16.8}{\rmdefault}{\mddefault}{\updefault}{ 0}%
}}}
\put(6751,-1861){\makebox(0,0)[lb]{\smash{\SetFigFont{17}{20.4}{\rmdefault}{\bfdefault}{\updefault}{ -3}%
}}}
\put(6226,-1261){\makebox(0,0)[lb]{\smash{\SetFigFont{17}{20.4}{\rmdefault}{\bfdefault}{\updefault}{ -1}%
}}}
\put(6076,314){\makebox(0,0)[lb]{\smash{\SetFigFont{17}{20.4}{\rmdefault}{\mddefault}{\updefault}{ C  }%
}}}
\put(3451,-1486){\makebox(0,0)[lb]{\smash{\SetFigFont{17}{20.4}{\rmdefault}{\mddefault}{\updefault}{ R}%
}}}
\put(8926,-1861){\makebox(0,0)[lb]{\smash{\SetFigFont{17}{20.4}{\rmdefault}{\mddefault}{\updefault}{ X}%
}}}
\put(9151,-1936){\makebox(0,0)[lb]{\smash{\SetFigFont{14}{16.8}{\rmdefault}{\mddefault}{\updefault}{ 0}%
}}}
\put(6601,-2536){\makebox(0,0)[lb]{\smash{\SetFigFont{17}{20.4}{\rmdefault}{\bfdefault}{\updefault}{ -2}%
}}}
\put(6376,314){\makebox(0,0)[lb]{\smash{\SetFigFont{17}{20.4}{\rmdefault}{\mddefault}{\updefault}{ (-6)}%
}}}
\put(8626,-1186){\makebox(0,0)[lb]{\smash{\SetFigFont{17}{20.4}{\rmdefault}{\mddefault}{\updefault}{     (-4)}%
}}}
\put(7576,-4186){\makebox(0,0)[lb]{\smash{\SetFigFont{17}{20.4}{\rmdefault}{\mddefault}{\updefault}{ X}%
}}}
\put(7801,-4336){\makebox(0,0)[lb]{\smash{\SetFigFont{14}{16.8}{\rmdefault}{\bfdefault}{\updefault}{ 2}%
}}}
\put(8026,-4186){\makebox(0,0)[lb]{\smash{\SetFigFont{17}{20.4}{\rmdefault}{\mddefault}{\updefault}{ (-4)}%
}}}
\put(6076,-5086){\makebox(0,0)[lb]{\smash{\SetFigFont{14}{16.8}{\rmdefault}{\mddefault}{\updefault}{ Figure 2}%
}}}
\put(5101,-4186){\makebox(0,0)[lb]{\smash{\SetFigFont{17}{20.4}{\rmdefault}{\mddefault}{\updefault}{ (-4)}%
}}}
\end{picture}

\par \vskip 1pc
The aim of this section is to characterize the alternating group
$A_6$ as the subgroup of $\text{Aut}(\mathcal D)$ fixing the six
vertices of a Coxeter-Dynkin diagram of type $A_{2}^{\oplus
2}\oplus A_{1}^{\oplus 2}$.

\begin{theorem}\label{theorem:leech}
\begin{list}{}{
\setlength{\leftmargin}{10pt}
\setlength{\labelwidth}{6pt}
}
\item[(1)] Up to $\text{Aut}(\mathcal D)$,
the graph $\mathcal R$ is the unique Coxeter-Dynkin diagram of
type $A_{2}^{\oplus 2}\oplus A_{1}^{\oplus 2}$, i.e. if $\mathcal R'
\subset \Pi_{2}$ is another set of $6$ Leech roots forming $A_{2}^{\oplus 2}\oplus A_{1}^{\oplus 2}$, then there is $\varphi \in\, \text{Aut}\,(\mathcal D)$
such that $\varphi(\mathcal R') = \mathcal R$.
\item[(2)] The group $\text{Aut}(\mathcal D, \mathcal R)$, the pointwise stabilizer of the set $\mathcal R$, is isomorphic to the alternating group
$A_{6}$.
\end{list}
\end{theorem}
\begin{proof} As before, we denote by $X$ the element of $\Lambda$ correspondng to $x \in \Pi_{2}$.
Since $Z = 0$, the stabilizer group of $z$ is contained in $\cdot 0$ under the identification
$\text{Aut}(\mathcal D)
= \cdot \infty$. The four vectors
$C- X_{0}$, $C- R_{0}$, $R_{0}$ and $X_{1}$  in $\Lambda$ form
the following diagram:

\par \vskip 3pc


\setlength{\unitlength}{3947sp}%
\begingroup\makeatletter\ifx\SetFigFont\undefined%
\gdef\SetFigFont#1#2#3#4#5{%
  \reset@font\fontsize{#1}{#2pt}%
  \fontfamily{#3}\fontseries{#4}\fontshape{#5}%
  \selectfont}%
\fi\endgroup%
\begin{picture}(5025,4903)(2101,-4694)
{ \thinlines
\put(2701,-211){\circle{150}}
}%
{ \put(4801,-3511){\circle{150}}
}%
{ \put(6976,-211){\circle{150}}
}%
{ \put(4801,-1486){\circle{150}}
}%
{ \put(2776,-211){\line( 1, 0){4050}}
}%
{ \put(2701,-361){\line( 2,-3){2019.231}}
}%
{ \put(6976,-361){\line(-2,-3){2042.308}}
}%
{ \put(2776,-286){\line( 5,-3){1875}}
}%
{ \put(6901,-286){\line(-5,-3){1930.147}}
}%
{ \put(4801,-1636){\line( 0,-1){1725}}
}%
\put(2101, 14){\makebox(0,0)[lb]{\smash{\SetFigFont{17}{20.4}{\rmdefault}{\mddefault}{\updefault}{ C}%
}}}
\put(2401, 14){\makebox(0,0)[lb]{\smash{\SetFigFont{17}{20.4}{\rmdefault}{\bfdefault}{\updefault}{ -}%
}}}
\put(2701, 14){\makebox(0,0)[lb]{\smash{\SetFigFont{17}{20.4}{\rmdefault}{\mddefault}{\updefault}{ X}%
}}}
\put(2926,-61){\makebox(0,0)[lb]{\smash{\SetFigFont{14}{16.8}{\rmdefault}{\mddefault}{\updefault}{ 0}%
}}}
\put(6901, 14){\makebox(0,0)[lb]{\smash{\SetFigFont{17}{20.4}{\rmdefault}{\mddefault}{\updefault}{ X}%
}}}
\put(7126,-61){\makebox(0,0)[lb]{\smash{\SetFigFont{14}{16.8}{\rmdefault}{\mddefault}{\updefault}{ 1}%
}}}
\put(4126,-1636){\makebox(0,0)[lb]{\smash{\SetFigFont{17}{20.4}{\rmdefault}{\mddefault}{\updefault}{ (-4)}%
}}}
\put(3226,-2011){\makebox(0,0)[lb]{\smash{\SetFigFont{17}{20.4}{\rmdefault}{\bfdefault}{\updefault}{ -2}%
}}}
\put(4576,-136){\makebox(0,0)[lb]{\smash{\SetFigFont{17}{20.4}{\rmdefault}{\bfdefault}{\updefault}{ -1}%
}}}
\put(5401,-886){\makebox(0,0)[lb]{\smash{\SetFigFont{17}{20.4}{\rmdefault}{\bfdefault}{\updefault}{ -1}%
}}}
\put(3526,-736){\makebox(0,0)[lb]{\smash{\SetFigFont{17}{20.4}{\rmdefault}{\bfdefault}{\updefault}{ -1}%
}}}
\put(4801,-2536){\makebox(0,0)[lb]{\smash{\SetFigFont{17}{20.4}{\rmdefault}{\bfdefault}{\updefault}{ 1}%
}}}
\put(6901,-661){\makebox(0,0)[lb]{\smash{\SetFigFont{17}{20.4}{\rmdefault}{\mddefault}{\updefault}{ (-4)}%
}}}
\put(2101,-511){\makebox(0,0)[lb]{\smash{\SetFigFont{17}{20.4}{\rmdefault}{\mddefault}{\updefault}{ (-4)}%
}}}
\put(4276,-1186){\makebox(0,0)[lb]{\smash{\SetFigFont{17}{20.4}{\rmdefault}{\mddefault}{\updefault}{ C}%
}}}
\put(4501,-1186){\makebox(0,0)[lb]{\smash{\SetFigFont{17}{20.4}{\rmdefault}{\bfdefault}{\updefault}{ -}%
}}}
\put(4726,-1186){\makebox(0,0)[lb]{\smash{\SetFigFont{17}{20.4}{\rmdefault}{\mddefault}{\updefault}{ R}%
}}}
\put(4951,-1261){\makebox(0,0)[lb]{\smash{\SetFigFont{14}{16.8}{\rmdefault}{\mddefault}{\updefault}{ 0}%
}}}
\put(4876,-3961){\makebox(0,0)[lb]{\smash{\SetFigFont{17}{20.4}{\rmdefault}{\mddefault}{\updefault}{ R}%
}}}
\put(5101,-4036){\makebox(0,0)[lb]{\smash{\SetFigFont{14}{16.8}{\rmdefault}{\mddefault}{\updefault}{ 0}%
}}}
\put(4351,-3961){\makebox(0,0)[lb]{\smash{\SetFigFont{17}{20.4}{\rmdefault}{\mddefault}{\updefault}{ (-4)}%
}}}
\put(4426,-4636){\makebox(0,0)[lb]{\smash{\SetFigFont{14}{16.8}{\rmdefault}{\mddefault}{\updefault}{ Figure 3}%
}}}
\put(6076,-1936){\makebox(0,0)[lb]{\smash{\SetFigFont{17}{20.4}{\rmdefault}{\bfdefault}{\updefault}{ -2}%
}}}
\end{picture}

\vskip 1pc \noindent
The four vectors generate a 4-dimensional lattice
containing exactly nine $(-4)$-vectors and six $(-6)$-vectors,
yielding a 4-dimensional $S$-lattice of type $2^93^6$ given in
[Cu, Page 554]. By [ibid, main theorem] (see
also Page 567), such a diagram is unique up to $\cdot 0$.

Let us consider the subgroup of $\cdot 0$ pointwise stabilizing $\{C- X_{0},
C- R_{0}, R_{0}, X_{1}\}$. This subgroup is obviously
the same as the one pointwise stabilizing $\{C, X_{0}, R_{0}, X_{1}\}$.
We denote this subgroup by $G_0$.

Note that $C^2 = -6$. So $G_0$ is a subgroup of $\cdot 3 < \cdot
0$.  By [ibid], $G_{0} = 3^4:A_{6}$, a semi-direct product of the
elementary abelian group $3^4 (= C_{3}^{\oplus 4})$ and $A_6$. By
[ibid], the lattice $M$ generated by $\{C, X_{0}, R_{0}, X_{1}\}$
is the maximal sublattice on which $G_{0}$ acts trivially.

Consider the set
$$\mathcal S:=\{V \in \Lambda \vert (V, C) = -3\,, (V, X_{0}) = (V, R_{0}) = (V, X_{1}) = -2\,,  V^{2} = -4\}.$$
Note that $X_2\in \mathcal S$. By calculating the intersection matrix,
we also see that $V \not\in \mathbf Q \langle C\, ,\, X_{0}\, ,\, R_{0}\, ,\, X_{1} \rangle$ for any $V \in \mathcal S$. In particular, $\mathcal S^{3^{4}:A_{6}} = \emptyset$ by the maximality of $M$
(this fact will be used in Lemma 2.5).

\begin{lemma}\label{lemma:leech} We have $\vert \mathcal S\vert = 81$. In other words, there are exactly $81$
Leech roots orthogonal to the five roots $c, z, x_0, r_0$ and $x_1$.
\end{lemma}
\begin{proof}
Since $V^{2} = -4$, we have $V=((\pm 2)^{8}, 0^{16})$, $(\mp 3, (\pm 1)^{23})$, or $((\pm 4)^2, 0^{22})$.

Case 1: $V=((\pm 2)^{8}, 0^{16})$. In this case $V=2\nu_K$, where $K$ is an octad with $K\ni \infty\,, 0$, $\vert K\cap K_0\vert= \vert K\cap K_1\vert= 4$. The number of such octads is 30.

Case 2: $V=(\mp 3, (\pm 1)^{23})$. In this case $V=4\nu_{\infty} +\nu_{\Omega}-2\nu_K$, where
$K=\{\infty, j_1, j_2,\dots, j_7\}$ is an octad with $K\not\ni 0$, $\vert K\cap K_0\vert= \vert K\cap K_1\vert= 2$. The number of such octads is 48.

Case 3: $V=((\pm 4)^{2}, 0^{22})$. In this case $V=4\nu_{\infty} +4\nu_{j}$, where $j\in \{1, 2, 3\}$. The number of such vectors is 3.
\end{proof}

Now we return to the proof of Theorem 2.3. The subgroup $G_{0}$ acts on the set $\mathcal S$. Moreover, $G_0$ acts on the following set $\mathcal S'$ of unordered pairs of vectors of norm $-4$:
$$\mathcal S' := \{\{W, C - W\} \vert W\in \Lambda\,,\, (W, C) = -3\,,\, W^{2} = -4 \}.$$
This action was studied by L. Finkelstein[Fi]. Note that $\vert
\mathcal S'\vert = 276$. By [ibid, Lemma 7.2], there is a subgroup
of $G_{0}$, isomorphic to $A_{6}$, whose action on $\mathcal S'$
has orbit decomposition $[1^{4}, 10, 20^{2}, 30^{2}, 36^{2},
45^{2}]$. By the fact that $A_{6}$ is simple, it also follows that
$G_{0} = 3^4:A_{6}$ for any $A_{6} < G_{0}$. Note also that
$\mathcal S \subset \mathcal S'$ in a natural manner :
$$\mathcal S\ni V\mapsto \{V, C-V\}\in\mathcal S'\,\, .$$
We also recall from [ibid, Lemma 4.8, Lemma 5.11] that our $3^{4}$
in $G_0$ is  a subgroup of the subgroup $3^{5}$ of $\cdot 3$,
which is unique up to conjugate, and that the action of $3^{5}$ on
$\mathcal S'$ has orbit decomposition $[3^{11}, 243]$.
\begin{lemma}\label{lemma:Tr} The action of $3^4$ on $\mathcal S$
is equivalent to the regular representation, i.e. the left action of $3^{4}$ on $3^{4}$.
\end{lemma}

\begin{proof} Since $243 = \vert 3^{5} \vert$,
the action of $3^{5}$ on the last orbit is equivalent to the
regular representation.  Therefore, the action of our $3^{4}$ on
$\mathcal S'$ has at most $33$ fixed points. In particular, its
action on $\mathcal S$ is non-trivial.

Since $\vert \mathcal S \vert = \vert 3^{4} \vert$, it suffices to show
that the action is transitive.

Let $\mathcal S = \cup_{k = 1}^{n} \mathcal O_{k}$ be the orbit
decomposition. Assume to the contrary that $n \geq 2$. Then $\vert
\mathcal O_{k} \mathcal \vert = 3^{m_{k}}$ for some integer $0
\leq m_{k} \leq 3$  and $\sum_{k = 1}^{n} 3^{m_{k}} = 81$. We may
assume that $m_{1} \geq m_{2} \geq \cdots \geq m_{n}$. Then
$m_{1}$ is either $1$, $2$, or $3$, i.e. $\vert \mathcal
O_{1} \vert$ is either $3$, $9$, or $27$.

Note that $A_{6}$ permutes these $n$ orbits $\mathcal O_k$. In what follows, we may often use
the {\it fact} that $A_{6}$ has no non-trivial homomorphism to $S_{m}$ with $m \leq 5$.
This {\it fact} is true because $A_{6}$ is simple.

First consider the case where $\vert \mathcal O_{1} \vert = 27$.
There are at most $3$ orbits of cardinality $27$. So, the orbit
$\mathcal O_{1}$ is stable under $A_{6}$. Thus, $\mathcal O_{1}$
must be a union of some orbits of the action of $A_{6}$ on
$\mathcal S'$. However, the shape of $[1^{4}, 10, 20^{2}, 30^{2},
36^{2}, 45^{2}]$ does not allow such a union, a contradiction.

Next consider the case where $\vert \mathcal O_{1} \vert = 9$.
There are at most $9$ orbits of cardinality $9$. By the above
{\it fact}, there are two cases: either there are exactly $m$ such
orbits, where $m$ is either $6$, $8$, $9$,  and $A_{6}$ acts on
the set $\{\mathcal O_{k} \}_{k =
1}^{m}$ of the $m$ orbits transitively, or at least one $\mathcal
O_{k}$ with $\vert \mathcal O_{k} \vert = 9$ is $A_{6}$-stable. In
the first case, if $m$ is $8$ or $9$,
then the stabilizer group of the action would be of order $360/8
= 45$ or $360/9 = 40$. However, $A_{6}$ has no subgroup of order $40$
or $45$ (see for instance the table of maximal subgroups of $A_{6}$), a
contradiction. Assume that $m = 6$ in the first case. Note that $3$
elements of $\mathcal S'$ corresponding to $X_{0}$,
$X_{1}$, $R_{0}$ are all fixed by $A_{6}$. So, there is in
$\mathcal S$ at most one element being fixed by $A_{6}$.
Thus $9 \times 6 = 54$ must be a sum of the entries of
$[1, 10, 20^{2}, 30^{2}, 36^{2}, 45^{2}]$. However, from the shape,
we see that there is no such sum, a contradiction.

In the second case, this $\mathcal O_{k}$ must be a
union of some orbits of the action of $A_{6}$ on $\mathcal S'$.
However, the shape of $[1^{4}, 10, 20^{2}, 30^{2}, 36^{2},
45^{2}]$ does not allow such a union, a contradiction.

Finally, consider the case where $\vert \mathcal O_{1} \vert = 3$.
As it is observed above, there is in
$\mathcal S$ at most one element being fixed by $A_{6}$. Note also
that there are only three ways to decompose $81$ into a sum of the
entries of $[1, 10, 20^{2}, 30^{2}, 36^{2}, 45^{2}]$. Namely, $81 = 1
+ 10 + 20 + 20 + 30$, $81 = 1 + 20 + 30 + 30$ and $81 = 36 + 45$.
In the first two cases, there is an element of $\mathcal S$ which
is fixed by $3^{4}:A_{6}$, a contradiction to the observation
immediately preceding Lemma 2.2. In the last case, all the orbits
$\mathcal O_{k}$ ($1 \le k \le 27$) are of order $3$, because $33 < 36$ and
$33 < 45$, and the orbit decomposition type of $A_{6}$ on
$\{\mathcal O_{k}\}_{k=1}^{27}$
is $[12, 15]$. Then, $A_{6}$ would have a subgroup $H$ of order
$360/12 = 30$. By the table of maximal subgroups of $A_{6}$, this
$H$ would then be an index $2$ subgroup of $A_{5}$, whence normal,
a contradiction to the fact that $A_{5}$ is simple.
\end{proof}

Thus, the configuration of pentagon formed by $\{C, X_{0}, R_{0},
X_{1}, X_{2}\}$ in $\Lambda$ is unique up to $\cdot 0$. This
implies the first assertion of Theorem 2.1. The second assertion
is now also clear, because the action of $3^4$ on $\mathcal S$ is
equivalent to the regular representation.
\end{proof}

The next proposition is also important in our proof of Theorems 3.1 and 5.1.

\begin{proposition}\label{proposition:SymR} Let $Sym(\mathcal R)$ be the group of symmetries of the Coxeter-Dynkin diagram of $\mathcal R$. Let $R\, (\simeq A_{2}^{\oplus 2}\oplus A_{1}^{\oplus 2})$ be the sublattice of $\Pi$ generated
by $\mathcal R$, the six Leech roots $c, z, x_{0}, r_{0}, x_{1}, x_{2}$.
Then
\begin{list}{}{
\setlength{\leftmargin}{10pt}
\setlength{\labelwidth}{6pt}
}
\item[(1)]
The natural homomorphism $Sym(\mathcal R)\longrightarrow O(A_R,
q_R)$ is an isomorphism, and $$Sym(\mathcal R)\simeq O(A_R,
q_R)\simeq D_8\times \mathbf Z/2,$$ where $D_8$ is the dihedral
group of order $8$.
\item[(2)]
Each element of $Sym(\mathcal R)$ is induced by an element in $\text{Aut}(\mathcal D)$.
\end{list}
\end{proposition}
\begin{proof}
The assertion (1) is obvious. Let us prove the assertion (2). One can find three Leech roots $u_1, u_2, u_3$ such that $x_2, u_1, c, z, u_2, r_{0}, x_{0}, u_3, x_{1}$ form a Coxeter-Dynkin diagram of type $A_{9}$ (Figure 4). For example,
$u_1=(4\nu_{0} +\nu_{\Omega}-2\nu_K, 1, 1)$, where
$K=\{0, 5, 12, 13, 16, 20, 21, 22\}$, $u_2=(4\nu_{0} +\nu_{\Omega}, 1, 2)$, $u_3=(\nu_{\Omega}-4\nu_{5}, 1, 1)$. There are also three Leech roots $u_2, v_1, v_2$ such that $c, z, u_2, r_{0}, x_{0}, v_1, v_2, x_{1}, x_2$ form one of type $D_{9}$ (Figure 5). For example, $v_1=(\nu_{\Omega}-4\nu_{7}, 1, 1)$, $v_2=(2\nu_K, 1, 1)$, where
$K=\{\infty, 0, 6, 7, 10, 12, 15, 18\}$.

\par \vskip 3pc

\setlength{\unitlength}{3947sp}%
\begingroup\makeatletter\ifx\SetFigFont\undefined%
\gdef\SetFigFont#1#2#3#4#5{%
  \reset@font\fontsize{#1}{#2pt}%
  \fontfamily{#3}\fontseries{#4}\fontshape{#5}%
  \selectfont}%
\fi\endgroup%
\begin{picture}(5625,1718)(751,-1394)
\thinlines
{ \put(6226, 89){\vector( 0,-1){0}}
\put(3639, 89){\oval(5174,450)[tr]}
\put(3639, 89){\oval(5176,450)[tl]}
\put(1051, 89){\vector( 0,-1){0}}
}%
{ \put(976,-136){\circle{108}}
}%
{ \put(2326,-136){\circle{108}}
}%
{ \put(3001,-136){\circle{108}}
}%
{ \put(4351,-136){\circle{108}}
}%
{ \put(5026,-136){\circle{108}}
}%
{ \put(6301,-136){\circle{108}}
}%
{ \put(1726,-136){\line( 1, 0){525}}
}%
{ \put(2401,-136){\line( 1, 0){525}}
}%
{ \put(3751,-136){\line( 1, 0){525}}
}%
{ \put(4426,-136){\line( 1, 0){525}}
}%
{ \put(5701,-136){\line( 1, 0){525}}
}%
{ \put(1051,-136){\line( 1, 0){525}}
}%
{ \put(5101,-136){\line( 1, 0){525}}
}%
{ \put(3076,-136){\line( 1, 0){525}}
}%
\put(1576,-286){\makebox(0,0)[lb]{\smash{\SetFigFont{14}{16.8}{\rmdefault}{\mddefault}{\updefault}{ *}%
}}}
\put(3601,-286){\makebox(0,0)[lb]{\smash{\SetFigFont{14}{16.8}{\rmdefault}{\mddefault}{\updefault}{ *}%
}}}
\put(5626,-286){\makebox(0,0)[lb]{\smash{\SetFigFont{14}{16.8}{\rmdefault}{\mddefault}{\updefault}{ *}%
}}}
\put(751,-436){\makebox(0,0)[lb]{\smash{\SetFigFont{17}{20.4}{\rmdefault}{\mddefault}{\itdefault}{ x}%
}}}
\put(1426,-436){\makebox(0,0)[lb]{\smash{\SetFigFont{17}{20.4}{\rmdefault}{\mddefault}{\itdefault}{ u}%
}}}
\put(3451,-436){\makebox(0,0)[lb]{\smash{\SetFigFont{17}{20.4}{\rmdefault}{\mddefault}{\itdefault}{ u}%
}}}
\put(4876,-436){\makebox(0,0)[lb]{\smash{\SetFigFont{17}{20.4}{\rmdefault}{\mddefault}{\itdefault}{ x}%
}}}
\put(5551,-436){\makebox(0,0)[lb]{\smash{\SetFigFont{17}{20.4}{\rmdefault}{\mddefault}{\itdefault}{ u}%
}}}
\put(6151,-436){\makebox(0,0)[lb]{\smash{\SetFigFont{17}{20.4}{\rmdefault}{\mddefault}{\itdefault}{ x}%
}}}
\put(3151,-1336){\makebox(0,0)[lb]{\smash{\SetFigFont{14}{16.8}{\rmdefault}{\mddefault}{\updefault}{ Figure 4}%
}}}
\put(1651,-586){\makebox(0,0)[lb]{\smash{\SetFigFont{14}{16.8}{\rmdefault}{\mddefault}{\updefault}{ 1}%
}}}
\put(976,-586){\makebox(0,0)[lb]{\smash{\SetFigFont{14}{16.8}{\rmdefault}{\mddefault}{\updefault}{ 2}%
}}}
\put(3676,-586){\makebox(0,0)[lb]{\smash{\SetFigFont{14}{16.8}{\rmdefault}{\mddefault}{\updefault}{ 2}%
}}}
\put(5101,-586){\makebox(0,0)[lb]{\smash{\SetFigFont{14}{16.8}{\rmdefault}{\mddefault}{\updefault}{ 0}%
}}}
\put(5776,-586){\makebox(0,0)[lb]{\smash{\SetFigFont{14}{16.8}{\rmdefault}{\mddefault}{\updefault}{ 3}%
}}}
\put(6376,-586){\makebox(0,0)[lb]{\smash{\SetFigFont{14}{16.8}{\rmdefault}{\mddefault}{\updefault}{ 1}%
}}}
\put(2101,-436){\makebox(0,0)[lb]{\smash{\SetFigFont{17}{20.4}{\rmdefault}{\mddefault}{\itdefault}{ c}%
}}}
\put(2776,-436){\makebox(0,0)[lb]{\smash{\SetFigFont{17}{20.4}{\rmdefault}{\mddefault}{\itdefault}{ z}%
}}}
\put(4126,-436){\makebox(0,0)[lb]{\smash{\SetFigFont{17}{20.4}{\rmdefault}{\mddefault}{\itdefault}{ r}%
}}}
\put(4276,-586){\makebox(0,0)[lb]{\smash{\SetFigFont{14}{16.8}{\rmdefault}{\mddefault}{\updefault}{ 0}%
}}}
\end{picture}

\par \vskip 1pc

\setlength{\unitlength}{3947sp}%
\begingroup\makeatletter\ifx\SetFigFont\undefined%
\gdef\SetFigFont#1#2#3#4#5{%
  \reset@font\fontsize{#1}{#2pt}%
  \fontfamily{#3}\fontseries{#4}\fontshape{#5}%
  \selectfont}%
\fi\endgroup%
\begin{picture}(5559,1769)(901,-1319)
\thinlines
{ \put(6076,-586){\vector(-1, 0){0}}
\put(6076,-136){\oval(750,900)[br]}
\put(6076,-136){\oval(750,900)[tr]}
\put(6076,314){\vector(-1, 0){0}}
}%
{ \put(3751,-136){\circle{108}}
}%
{ \put(3076,-136){\circle{108}}
}%
{ \put(1801,-136){\circle{108}}
}%
{ \put(1126,-136){\circle{108}}
}%
{ \put(5851,389){\circle{108}}
}%
{ \put(5888,-623){\circle{108}}
}%
{ \put(3151,-136){\line( 1, 0){525}}
}%
{ \put(1201,-136){\line( 1, 0){525}}
}%
{ \put(2476,-136){\line( 1, 0){525}}
}%
{ \put(1876,-136){\line( 1, 0){525}}
}%
{ \put(3826,-136){\line( 1, 0){525}}
}%
{ \put(4426,-136){\line( 1, 0){525}}
}%
{ \put(5701,389){\line(-3,-2){675}}
}%
{ \put(5071,-196){\line( 2,-1){750}}
}%
\put(901,-436){\makebox(0,0)[lb]{\smash{\SetFigFont{17}{20.4}{\rmdefault}{\mddefault}{\itdefault}{ c}%
}}}
\put(3526,-436){\makebox(0,0)[lb]{\smash{\SetFigFont{17}{20.4}{\rmdefault}{\mddefault}{\itdefault}{ x}%
}}}
\put(2851,-436){\makebox(0,0)[lb]{\smash{\SetFigFont{17}{20.4}{\rmdefault}{\mddefault}{\itdefault}{ r}%
}}}
\put(2176,-436){\makebox(0,0)[lb]{\smash{\SetFigFont{17}{20.4}{\rmdefault}{\mddefault}{\itdefault}{ u}%
}}}
\put(2326,-286){\makebox(0,0)[lb]{\smash{\SetFigFont{14}{16.8}{\rmdefault}{\mddefault}{\updefault}{ *}%
}}}
\put(4876,-286){\makebox(0,0)[lb]{\smash{\SetFigFont{14}{16.8}{\rmdefault}{\mddefault}{\updefault}{ *}%
}}}
\put(4276,-286){\makebox(0,0)[lb]{\smash{\SetFigFont{14}{16.8}{\rmdefault}{\mddefault}{\updefault}{ *}%
}}}
\put(1576,-436){\makebox(0,0)[lb]{\smash{\SetFigFont{17}{20.4}{\rmdefault}{\mddefault}{\itdefault}{ z}%
}}}
\put(4126,-436){\makebox(0,0)[lb]{\smash{\SetFigFont{17}{20.4}{\rmdefault}{\mddefault}{\itdefault}{ v}%
}}}
\put(4726,-436){\makebox(0,0)[lb]{\smash{\SetFigFont{17}{20.4}{\rmdefault}{\mddefault}{\itdefault}{ v}%
}}}
\put(2401,-511){\makebox(0,0)[lb]{\smash{\SetFigFont{14}{16.8}{\rmdefault}{\mddefault}{\updefault}{ 2}%
}}}
\put(3001,-511){\makebox(0,0)[lb]{\smash{\SetFigFont{14}{16.8}{\rmdefault}{\mddefault}{\updefault}{ 0}%
}}}
\put(3751,-511){\makebox(0,0)[lb]{\smash{\SetFigFont{14}{16.8}{\rmdefault}{\mddefault}{\updefault}{ 0}%
}}}
\put(4951,-511){\makebox(0,0)[lb]{\smash{\SetFigFont{14}{16.8}{\rmdefault}{\mddefault}{\updefault}{ 2}%
}}}
\put(4351,-511){\makebox(0,0)[lb]{\smash{\SetFigFont{14}{16.8}{\rmdefault}{\mddefault}{\updefault}{ 1}%
}}}
\put(5476,-811){\makebox(0,0)[lb]{\smash{\SetFigFont{17}{20.4}{\rmdefault}{\mddefault}{\itdefault}{ x}%
}}}
\put(5701,-961){\makebox(0,0)[lb]{\smash{\SetFigFont{14}{16.8}{\rmdefault}{\mddefault}{\updefault}{ 2}%
}}}
\put(2926,-1261){\makebox(0,0)[lb]{\smash{\SetFigFont{14}{16.8}{\rmdefault}{\mddefault}{\updefault}{ Figure 5}%
}}}
\put(5626, 89){\makebox(0,0)[lb]{\smash{\SetFigFont{17}{20.4}{\rmdefault}{\mddefault}{\itdefault}{ x}%
}}}
\put(5851,-61){\makebox(0,0)[lb]{\smash{\SetFigFont{14}{16.8}{\rmdefault}{\mddefault}{\updefault}{ 1}%
}}}
\end{picture}

\par \vskip 1pc
Thus, by Borcherds [Bo, Lemma 9.6, Theorem 9.5], there are two
isometries in $\text{Aut}(\mathcal D)$ whose restrictions on $R$
give rise to two involutions in $Sym(\mathcal R)$, one switching
the pair of $A_2$'s as well as the pair of $A_1$'s, and the other
switching the pair of $A_1$'s only. We need one more involution.
From Theorem 2.1 and Lemma 2.5, we know that the pointwise
stabilizer group $\text{Aut}(\mathcal D, \mathcal R)$ of the
pentagon formed by $\{C, X_0, R_0, X_1, X_2\}$ is a subgroup of
$\cdot 3$ isomorphic to $A_6$, and that the action of
$\text{Aut}(\mathcal D, \mathcal R)$ on the set $\mathcal S'$ has
only 4 fixed points $\{X_0, C-X_0\}$, $\{R_0, C-R_0\}$, $\{X_1,
C-X_1\}$, $\{X_2, C-X_2\}$ (the maximality of the lattice $M$
preceding Lemma 2.4). By [Fi, Lemma 7.2], the normalizer $N_{\cdot
3}(\text{Aut}(\mathcal D, \mathcal R))$ in $\cdot 3$ is isomorphic
to $\mathbf Z/2 \times \text{Aut}(A_6)$, and hence has quotient
$N_{\cdot 3}(\text{Aut}(\mathcal D, \mathcal
R))/\text{Aut}(\mathcal D, \mathcal R) \simeq \mathbf Z/2 \times
(\mathbf Z/2)^2$. So, it is easy to see that the first factor
$\mathbf Z/2$ sends the pentagon to its dual pentagon $\{C, C-X_0,
C-R_0, C-X_1, C-X_2\}$, and the second $(\mathbf Z/2)^2$ gives the
full symmetry group of the pentagon. From this it is obvious to
see that there is an element of
 $\cdot 3$ interchanging $R_0$ and $X_0$, and leaving $C$, $X_1$,
 and $X_2$ fixed. Since $\cdot 3 < \text{Aut}(\mathcal D)$,
 this gives an isometry in $\text{Aut}(\mathcal D)$ whose restriction on
 ${\mathcal R}$ is the involution switching $r_{0}$ and $x_{0}$ while
 leaving the remaining four roots fixed. Clearly this involution, together with the previous two, generate $Sym(\mathcal R)$.
 This completes the proof.
\end{proof}

By the above lemma, one can choose an isometry $\phi_4\in \text{Aut}(\mathcal D)$
inducing the order 4 element of $Sym(\mathcal R)$:
$$c\mapsto x_0\mapsto z\mapsto r_0\mapsto c\,,\quad x_1\leftrightarrow x_2.$$
The choice of $\phi_4$ is unique up to $\text{Aut}(\mathcal D, \mathcal R)\simeq A_6$.
Also $\phi_4^4$ is in $A_6$.

We take a basis of the discriminant group $A_R\simeq (\mathbf
Z/6)^2$:
$$e_1=(c+2z)/3+(r_0+2x_0)/3+x_1/2\,,\quad e_2=(c+2z)/3-(r_0+2x_0)/3+x_2/2,$$
whose intersection matrix is
$$\begin{pmatrix}
1/6&0\cr
0&1/6
\end{pmatrix}\,\, .$$

Then $\phi_4$ acts on $A_R$ as follows:
$$\phi_4 : e_1\mapsto e_2\,,\quad e_2\mapsto -e_1.$$

\begin{definition} \label{definition:alt}
(Definition of $\tilde{A_{6}}$)
 The subgroup of $\text{Aut}(\mathcal D)$ generated
by $\phi_4$ and $\text{Aut}(\mathcal D, \mathcal R)$ is denoted by
$\tilde{A_{6}}$. This is a uniquely determined group.
\end{definition}

We denote by $A_{6}.\mu_{4}$ {\it any} group $G$ which falls into the exact sequence:

$$1 \longrightarrow A_{6} \longrightarrow G \longrightarrow \mu_{4} \longrightarrow 1\,\, $$

There are many isomorphism classes of groups of the form
$A_{6}.\mu_{4}$. The direct product group $A_{6} \times \mu_{4}$
is clearly one of $A_{6}.\mu_{4}$. Our group $\tilde{A}_{6}$ is
also an extension of $A_6$ by $\mu_{4}$, and is one of
$A_{6}.\mu_{4}$. Indeed, $\text{Ker}\, f = A_{6}$ and $\text{Im}\,
f  \simeq \mathbf Z/4$ for the natural homomorphism $f :
\tilde{A}_{6} \longrightarrow Sym(\mathcal R)$.

We close this section by introducing a very important vector $h
\in R^{\perp}_{\Pi}$, the orthogonal complement of $R$ in $\Pi$,
and its properties. This vector will play a crucial role in
Sections 3 and 5:

\begin{lemma} \label{lemma:isom}
Let $w'$ be the orthogonal projection of
the Weyl vector $w$ onto
$R^{\perp}_{\Pi} \otimes \mathbf Q$, or more explicitly, $w'$ is
the vector such that
$$w = w' + w_{R},\quad \text{where}\quad w_{R} = -(c +
z + r_{0} + x_{0} + (x_{1} + x_{2})/2).$$
We set $h := 2 w'$. Then

\begin{list}{}{
\setlength{\leftmargin}{10pt}
\setlength{\labelwidth}{6pt}
}
\item[(1)]
$h^2=20$, $h \in R^{\perp}_{\Pi}$ and $h$ is primitive in $R^{\perp}_{\Pi}$.
\item[(2)]
$\varphi(h) = h$ for any element $\varphi$ of $\text{Aut}(\mathcal D)$ such that $\varphi(\mathcal R) = \mathcal R$.
\item[(3)]
There is no element $x \in R^{\perp}_{\Pi}$ such that $x^{2} = -2$
and $(x, h) = 0$.
\end{list}
\end{lemma}

\begin{proof}
The assertions (1) and (2) are obvious. Let us show (3). Recall that the Leech lattice $\Lambda$ has no element of norm $-2$ and that
$z$ and $w$ generate the second direct summand of the (original)
orthogonal decomposition $\Pi = \Lambda \oplus U$. Since
$$h^{\perp}_{R^{\perp}_{\Pi}} \simeq \langle R\,,\,w
\rangle^{\perp}_{\Pi} \subset\Lambda,$$ every $x \in R^{\perp}_{\Pi}$ with
$(x, h)=0$ must have norm $x^2< -2$.
\end{proof}

\section{Existence of a $K3$ surface with an $\tilde{A_{6}}$-action}

The goal of this section is to construct a triplet $(F,
\tilde{A}_{6}, \rho_{F})$ consisting of a $K3$ surface $F$ and a
faithful $\tilde{A}_{6}$-action $\rho_{F} : \tilde{A}_{6} \times F
\longrightarrow F$ (Theorem (3.1)), and then to give an explicit
description of $F$ (Proposition (3.5)). In the next two sections, it
turns out that such a triplet $(F, \tilde{A}_{6}, \rho_{F})$ is
actually unique up to isomorphisms.

Before stating our main result of this section, we recall some facts about such triplets.

Throughout this note, by a $K3$ surface, we mean a
simply-connected compact complex surface $X$ admitting a nowhere
vanishing global holomorphic $2$-form $\omega_{X}$. The second
cohomology group $H^{2}(X, \mathbf Z)$ together with a cup product
becomes an even unimodular lattice of index $(3, 19)$ and is
isomorphic to the so-called $K3$ lattice $U^{\oplus 3} \oplus
E_{8}^{\oplus 2}$, where $E_8$ is the {\it negative} definite even
unimodular lattice of rank 8. We denote by $S(X)$ the
N\'eron-Severi lattice of $X$. This is a primitive sublattice of
$H^{2}(X, \mathbf Z)$ generated by the classes of line bundles. We
denote by $T(X)$ the transcendental lattice of $X$, i.e. the
minimal primitive sublattice whose $\mathbf C$-linear extension
contains the class $\omega_{X}$, or equivalently $T(X) =
S(X)^{\perp}$ in $H^{2}(X, \mathbf Z)$. If $X$ is projective, then
$S(X) \cap T(X) = \{0\}$ and $S(X) \oplus T(X)$ is a finite-index
sublattice of $H^{2}(X, \mathbf Z)$.

Let $(X, G, \rho)$ be a triplet consisting of a $K3$ surface, a
finite group $G$ and a faithful action $\rho : G \times X
\longrightarrow X$. Then $G$ has a $1$-dimensional representation
on $H^{0}(X, \Omega_{X}^{2}) = \mathbf C \omega_{X}$ given by
$g^{*}\omega_{X} = \alpha(g)\omega_{X}$, and we have an exact
sequence, called the basic sequence:
$$1 \longrightarrow G_{N} := \text{Ker}\, \alpha \longrightarrow G
\mapright{\alpha}
\mu_{I} \longrightarrow 1\,\,  .$$ We call $G_{N}$ the symplectic
part and $\mu_{I} := \langle \zeta_{I} \rangle$ (resp. $I$), where
$\zeta_{I} = \text{exp}(2 \pi \sqrt{-1}/I)$, the transcendental
part (resp. the transcendental value) of the action $\rho$. We
note that if $A_{6}.\mu_{4}$ acts faithfully on a $K3$ surface
then $G_{N} \simeq A_{6}$ and the transcendental part is
isomorphic to $\mu_{4}$. This follows from the fact that $A_{6}$
is simple and also maximal among all finite groups acting on a
$K3$ surface symplectically [Mu]. We also note that $X$ is
projective if $I \geq 2$ [Ni1].

We say that 2 triplets $(X, G, \rho)$ and $(X', G', \rho')$ are isomorphic if
there are a group isomorphism $f : G' \simeq G$ and an isomorphism
$\varphi : X' \simeq X$ such that the following diagram commutes:
$$
\begin{matrix}
G \times X & \mapright{\rho} & X \cr
\mapupleft{f \times \varphi \; } &  & \mapup{\varphi} \cr
G' \times X' & \mapright{\rho'} & X' \cr
\end{matrix}
$$

The aim of this section is to show the following:
\begin{theorem} \label{theorem:exist}
There is a triplet $(F, \tilde{A_{6}}, \rho_{F})$ consisting of a $K3$ surface $F$ and a faithful group action
$\rho_{F} : \tilde{A}_{6} \times F \longrightarrow F$
of $\tilde{A_{6}}$ on $F$. Here $\tilde{A}_{6}$ is the group
defined in (2.7).
\end{theorem}
\begin{proof}

Let $F$ be a $K3$ surface such that the transcendental lattice $T(F) = \mathbf Z \langle t_{1}, t_{2} \rangle$ has the intersection matrix
$$\begin{pmatrix}
6&0\cr
0&6
\end{pmatrix}\,\, ,$$
and $\omega_{F} := t_{1} + \sqrt{-1} t_{2}$ is a holomorphic $2$-form of $F$.
Such a $K3$ surface exists and is unique [SI]. We claim that the surface $F$ admits an action of $\tilde{A_{6}}$ as a group of automorphisms.

In order to find an action of $\tilde{A_{6}} = \langle \text{Aut}(\mathcal D, \mathcal R), \phi_{4} \rangle$ on $F$, we first relate the N\'eron-Severi lattice $S(F)$ to the lattices $\Pi = \Lambda \oplus U$ and $R$. Here and hereafter, we shall freely use the lattices and their elements introduced in Section 2.

The Picard lattice $S(F)$ is isometric to
$U\oplus E_8\oplus E_8\oplus \langle -6 \rangle \oplus
\langle -6 \rangle$.
Indeed, one has $(A_{S(F)}, q_{S(F)}) \simeq (A_{T(F)}, -q_{T(F)})$,
and the genus of $S(F)$ is the single element set $\{S(F)\}$
[Ni2, Theorem 1.14.2]. We set $L := H^{2}(F, \mathbf Z)$. Since
$(A_{S(F)}, q_{S(F)})\simeq (A_R, -q_R)$, one has also an isomorphism
$$\Phi : S(F)\simeq R^{\perp}_{\Pi} \subset \Pi\,\, ,$$
and the diagram (depending on $\Phi$):
$$T(F) \subset L \supset S(F) \simeq_{\Phi} R^{\perp}_{\Pi} \subset \Pi \supset R\,\, .$$

Since $L$ and $\Pi$ are both unimodular, these primitive inclusions and $\Phi$ naturally induce the isomorphisms, depending on $\Phi$, of the discriminant
groups:

$$(A_{T(F)}, q_{T(F)}) \simeq (A_{S(F)}, -q_{S(F)}) \simeq_{\Phi}
(A_{R^{\perp}_{\Pi}}, -q_{R^{\perp}_{\Pi}}) \simeq (A_{R}, q_{R})\,\, .$$

We also recall that $t_{1}/6$ and $t_{2}/6$ (resp.
$e_{1}$ and $e_{2}$ defined in the previous section) are generators of $A_{T(F)}$ (resp. $A_{R}$) with intersection form
$$\begin{pmatrix}
1/6&0\cr
0&1/6
\end{pmatrix}\,\, .$$

Next, we shall transfer the group action of $\tilde{A}_{6}$ on $\Pi$ to an effective Hodge isometric action on $L$. For this, it is more convenient to choose a special isomorphism $\Phi$ given by the next Lemma:

\begin{lemma} \label{lemma:isom}
There is an isomorphism $\Phi : S(F) \simeq R^{\perp}_{\Pi}$ such that:
\begin{list}{}{
\setlength{\leftmargin}{10pt}
\setlength{\labelwidth}{6pt}
}
\item[(1)]
$H := \Phi^{-1}(h)$ is an ample class of $F$, and
\item[(2)]
under the isomorphism $(A_{T(F)}, q_{T(F)}) \simeq (A_{R}, q_{R})$ above,
we have $t_{1}/6 \leftrightarrow e_{1}$ and $t_{2}/6 \leftrightarrow e_{2}$.
\end{list}
Here $h \in R^{\perp}_{\Pi}$ is
the vector defined in Lemma (2.8).
\end{lemma}

\begin{proof}

Let us choose any isomorphism $\Phi_{0} : S(F) \simeq R^{\perp}_{\Pi}$.
By Lemma (2.8), there is no vector $x \in S(F)$ such that $(x, \Phi_{0}^{-1}(h)) = 0$ and $(x^{2}) = -2$. This means that a product
$\sigma$ of $(-2)$-reflections of $S(F)$ and $-1_{S(F)}$
transforms $\Phi_{0}^{-1}(h)$ to an ample divisor class $H$ on $F$.
So, the new isomorphism $\Phi_{1} := \Phi_{0} \circ \sigma^{-1} : S(F) \simeq R^{\perp}_{\Pi}$ enjoys the property (1). Let $e_{1}'$ and $e_{2}'$ be the generators of $A_{R}$ corresponding to $t_{1}/6$ and $t_{2}/6$ under $\Phi_{1}$.
Then, by Proposition (2.6), there is $\eta \in \text{Aut}(\mathcal D)$
such that $\eta(\mathcal R) = \mathcal R$ and $\eta(e_{1}') = e_{1}$ and
$\eta(e_{2}') = e_{2}$ (on $A_{R}$). Note also that $\eta(h) = h$
by Lemma (2.8)(2).
Now, the new isomorphism
$\Phi := \eta \circ \Phi_{1} : S(F) \simeq R^{\perp}_{\Pi}$ satisfies
the properties (1) and (2).
\end{proof}

{\it Using this} $\Phi${\it, we shall identify $S(F) = R^{\perp}_{\Pi}$ and $H = h$ from now on}:

$$T(F) \subset L \supset S(F) = R^{\perp}_{\Pi} \subset \Pi \supset R\,\, .$$

Let us construct an action of $\tilde{A_{6}} = \langle \text{Aut}(\mathcal D, \mathcal R), \phi_{4}\rangle$ on $L$. First, observe that the group $\tilde{A_{6}}$
acts on $R^{\perp} = S(F)$, faithfully, i.e.
$\tilde{A_{6}}$ can be viewed as a subgroup of $O(S(F))$.
Indeed, if $\varphi \vert S(F) = id$ for $\varphi \in \tilde{A_{6}}$,
then $\varphi \vert {A_{R}} = id$ as well. Since $Sym(\mathcal R) \simeq O(A_{R}, q_{R})$, we have then $\varphi \vert R = id$ and consequently $\varphi = id$ on $\Pi$.

Let us define the isometry $\psi_{4} \in O(T(F))$
by $t_1\mapsto t_2$, $t_2\mapsto -t_1$. Then, by using Nikulin [Ni2,
Corollary 1.5.2], one can find an isometry $\tilde{\phi}_{4} \in O(L)$ such that $\tilde{\phi}_{4} \vert S(F) = \phi_{4} \vert S(F)$ and
$\tilde{\phi}_{4} \vert T(F) = \psi_{4}$.
Let $\phi \in \text{Aut}(\mathcal D, \mathcal R)$. Then, $\phi \vert A_{S(F)} = id$ and we have an isometry $\tilde{\phi} \in O(L)$ such that $\tilde{\phi} \vert S(F) = \phi \vert S(F)$ and $\tilde{\phi} \vert T(F) = id$.

Let $\tilde{A_{6}}'$ be the subgroup of $O(L)$ generated by $\tilde{\phi}_{4}$ and these $\tilde{\phi}$:

$$\tilde{A_{6}}' := \langle \tilde{\phi}\, (\forall \phi \in \text{Aut}(\mathcal D, \mathcal R))\,\, ,\,\, \tilde{\phi}_{4} \rangle < O(L)\,\, .$$

\begin{lemma} \label{lemma:isom}
\begin{list}{}{
\setlength{\leftmargin}{10pt}
\setlength{\labelwidth}{6pt}
}
\item[(1)]
The natural homomorphism $\iota : \tilde{A_{6}}' \longrightarrow
\tilde{A_{6}} (< O(S(F))$, induced by the restriction of the action on $S(F)$,
is an isomorphism.
\item[(2)]
Each element of $\tilde{A_{6}}'$ is an effective Hodge isometry of $L$.
\end{list}
\end{lemma}

\begin{proof} As we have already observed, $\phi \vert S(F)$
($\phi \in \text{Aut}(\mathcal D, \mathcal R)$) and $\phi_{4} \vert S(F)$
generate
$\tilde{A_{6}} (< O(S(F)))$. Thus, $\iota$ is surjective. Let $\varphi \in \tilde{A_{6}}'$ such that $\varphi \vert S(F) = id$. Then, $\varphi \vert A_{S(F)} = id$ and $\varphi \vert A_{T(F)} = id$ as well. Observe that the natural homomorphism $D_{8} \simeq O(T(F)) \longrightarrow O(A_{T(F)}, q_{T(F)}) \simeq D_{8} \times C_{2}$ is injective. Then, $\varphi \vert T(F) = id$ and hence
$\varphi \vert L = id$. This means that $\iota$ is also injective.

Let us show the assertion (2). It suffices to check it for the
generators. It is clear that $\tilde{\phi}$ ($\phi \in
\text{Aut}(\mathcal D, \mathcal R)$) preserves the Hodge
decomposition of $L$. By $\omega_{F} = t_{1} + \sqrt{-1}t_{2}$, we
have $\psi_{4}(\omega_{F}) = -\zeta_{4}\omega_{F}$. Thus
$\tilde{\phi}_{4}$ also preserves the Hodge decomposition of $L$.
In addition, our group $\tilde{A_{6}}' \simeq \tilde{A_{6}}$,
being a subgroup of $\text{Aut}(\mathcal D)$, fixes the Weyl
vector $w$ and the set $\mathcal R$ by the definition, whence it
fixes $H = h$ by Lemma (2.8). Since $H$ is ample, the action of
$\tilde{A_{6}}'$ on $L$ is also effective.
\end{proof}
Thus the group $\tilde{A_{6}}' \simeq \tilde{A_{6}}$
realizes as a group of automorphisms of $F$ by the Torelli Theorem for
$K3$ surfaces [PSS], [BR]. This completes the proof of Theorem (3.1).
\end{proof}

\begin{remark} \label{}
It will turn out that $S(F)^{\tilde{A}_{6}} =
S(F)^{A_{6}} = \mathbf Z H$ by Proposition (4.5) and Lemma (2.8).

\end{remark}

We shall close this section by giving an explicit equation of a canonical model of $F$ in $\mathbf P^{1} \times \mathbf P^{2}$.

\begin{proposition} \label{proposition:eqn}
Let $F$ be a $K3$ surface constructed in Theorem (3.1), i.e. the unique $K3$ surface whose transcendental lattice has the intersection matrix
$$\begin{pmatrix}
6&0\cr
0&6
\end{pmatrix}\,\, .$$
\begin{list}{}{
\setlength{\leftmargin}{10pt}
\setlength{\labelwidth}{6pt}
}
\item[(1)]
The $K3$ surface $F$ is isomorphic to the minimal resolution of
the double cover $\overline{F}$ of the (rational) elliptic modular surface $E$
with level $3$ structure. The double cover is branched along two
of a total of $4$ singular fibres of the same type $I_3$ and $\overline{F}$ has
$6$ ordinary double points.
\item[(2)]
$\overline{F}$ is isomorphic to a surface in $\mathbf P^{1} \times
\mathbf P^{2}$ given by the following equation, where
$([S:T],[X:Y:Z])$ are coordinates of $\mathbf P^{1} \times \mathbf
P^{2}$:
$$S^2(X^3+Y^3+Z^3) - 3(S^2+T^2)XYZ = 0.$$
\end{list}
\end{proposition}
\begin{proof}
Let $E := \{\lambda (X^3+Y^3+Z^3) - 3 \mu XYZ = 0\} \subset
\mathbf P^{1} \times \mathbf P^{2}$ be the (rational and smooth)
elliptic modular surface with level 3 structure. It is easy to see
that the elliptic fibration $E \longrightarrow \mathbf
P^{1}$, induced from the projection $\mathbf P^{1} \times \mathbf
P^{2} \rightarrow \mathbf P^1$, has exactly four singular fibres
of the same type $I_3$ lying over the points $[0:1]$,
$[1:\zeta_3]$, $[1:\zeta_3^2]$ and $[1:1]$. Let $\mathbf P^1
\rightarrow \mathbf P^1$, $[S:T] \mapsto [S^2: S^2+T^2]$, be the
double cover branched at $[0:1]$ and $[1:1]$. Then the pull back
${\overline F}$ of $E$ in the fibre product $\mathbf P^{1} \times
\mathbf P^{2} \cong \mathbf P^1 \times_{\mathbf P^1} \mathbf P^2$
is given by the equation:
$$S^2(X^3+Y^3+Z^3) - 3(S^2+T^2)XYZ = 0.$$
Now ${\overline F} \rightarrow E$ is branched along the fibres
$E_1$ and $E_{\infty}$. Here we let $E_t$ (resp. $E_{\infty}$) be
the fibre lying over $[1:t]$ (resp. $[0:1]$). Note that
${\overline F}$ has six singular points of Dynkin type $A_1$ lying
over the six intersection points in $E_1$ and $E_{\infty}$ and is
smooth everywhere else. Let $F \rightarrow {\overline F}$ be the
minimal resolution of these $6$ singular points. We shall show
that this $F$ is isomorphic to its namesake constructed in Theorem
(3.1).

\par
The adjunction formula shows that ${\overline F}$ has trivial
canonical line bundle. By a cohomology exact sequence and the
Kawamata-Viehweg vanishing theorem, we see that the irregularity
$q({\overline F}) = 0$. So $F$ is a $K3$ surface. The elliptic
fibration on $E$ lifts to an elliptic fibration $F \rightarrow
\mathbf P^1$ with a section. Now $E_1$ and $E_{\infty}$ lift to
singular fibres of type $I_6$, while each of the singular fibres
$E_{\zeta^i}$ ($i = 1,2$) splits into two singular fibres of the
same type $I_3$. Thus, by [Sh], $\rho(X) \geq 2 + 4 \cdot 2 + 2
\cdot 5 = 20$, whence $\rho(X) = 20$. According to [SZ, Table 2,
No. 5], the transcendental lattice $T(F)$ has the intersection
matrix $(a_{ij})$ of rank 2 with $a_{11}=a_{22} = 6$ and $a_{12} =
a_{21} = 0$. So the surface $F$ is exactly the same $K3$ surface
as in Theorem (3.1) by [SI]. This proves the proposition.
\end{proof}

\begin{remark} \label{}
(1) The surface $F$ is not a Kummer surface, as its transcendental
lattice is not the double of an even lattice [Mor].

(2) The surface $F$ is the universal double cover of an Enriques
surface, i.e. has a fixed point free involution. This can be seen
indirectly by the criterion in [Ke].

\end{remark}

\section{Uniqueness of the $K3$ surface admitting
an $A_{6}.\mu_{4}$-action}

In this section, we shall show the uniqueness of the $K3$ surface admitting an $A_{6}.\mu_{4}$-action and the maximality of the extension of
$A_{6}$ by $\mu_{4}$. Our main results of this section are Propositions (4.1) and (4.5).

In what follows, we set $L := H^{2}(X, \mathbf Z)$ for a $K3$
surface $X$. We define
$$L^{G_{N}} :=\{x\in L\,|\,g^*x=x\,\,\text{for all}\,\, g\in G_{N}\}$$
$$L_{G_{N}} :=(L^{G_{N}})^{\perp}_L =\{x\in L\,|\, (x,y)=0\,\,\text{for all}\,\, y\in L^{G_{N}}\}.$$

\begin{proposition} \label{proposition:trans}
Let $G$ be a finite group acting faithfully on a $K3$ surface $X$. Assume that
$A_{6} < G$ and $I \geq 2$, where $I$ is the transcendental value. Then we have
\begin{list}{}{
\setlength{\leftmargin}{10pt}
\setlength{\labelwidth}{6pt}
}
\item[(1)]
$G_{N} = A_{6}$.
\item[(2)]
$\text{rank}\, L^{G_{N}} = 3$. In particular, $S(X)^{G_{N}} =
\mathbf Z H$, where $H$ is an ample class, and $\text{rank}\, T(X)
= 2$.
\item[(3)]
$I = 2$, or $4$.
\end{list}
\end{proposition}
\begin{proof}
As we remarked in Section 3, the statement (1) follows from the
fact that $[A_6, A_6] = A_6$ and the maximality of $A_6$ as a
symplectic $K3$ group [Mu]. Since $X$ is projective by $I \geq 2$,
the second statement of (2) follows from the first one. Let us
prove the first statement of (2).

Recall that the order structure of $A_{6}$ is as follows:
$$
\begin{tabular}{*{8}{|c}{|}}
\hline
\text{order[conjugacy class]} & 1 [1A] & 2 [2A] & 3 [3A] & 3 [3B] & 4 [4A]
& 5 [5A] & 5 [5B]\\
\hline
\text{cardinality}  & 1 & 45 & 40 & 40 & 90 & 72 & 72\\
\hline
\end{tabular}
$$
Moreover, by [Ni1], the number of the fixed points of the symplectic
action is as follows.
$$
\begin{tabular}{*{9}{|c}{|}}
\hline
$ord(g)$ & 1 & 2 & 3 & 4 & 5 & 6 & 7 & 8 \\
\hline
$\vert X^{g} \vert$ & $X$ & 8 & 6 & 4 & 4 & 2 & 3 & 2 \\
\hline
\end{tabular}
$$
Set $\tilde{H}(X, \mathbf Z) = H^{0}(X, \mathbf Z) \oplus
H^{2}(X, \mathbf Z) \oplus H^{4}(X, \mathbf Z)$.
Now, by applying the topological Lefschetz fixed point formula for $G_{N} = A_{6}$,
we calculate that
$$\text{rank}\, \tilde{H}(X, \mathbf Z)^{A_{6}}  =
\frac{1}{\vert A_{6} \vert}
\sum_{g \in A_{6}} \text{tr}(g^{*} \vert \tilde{H}(X, \mathbf Z))$$
$$ =
\frac{1}{360}(24 + 8\cdot 45 + 6 \cdot 80 + 4 \cdot 90 + 4 \cdot 144 )
= 5\,\, .$$
This implies the result.

Let us show the assertion (3). By (2), we have $\text{rank}\, T(X)
= 2$. Thus, $I = 2, 4, 3$ or $6$ because the Euler function
$\varphi(I)$ divides rank $T(X)$. It suffices to rule out the case
$I = 3$.

Let us first determine the irreducible decomposition of $S(X)
\otimes \mathbf C$ as $A_{6}$-modules. In the description, we use
Atlas notation for irreducible characters of $A_{6}$ in the Table
below. We also use the same letters for the representations.

\par \vskip 3pc


\begin{tabular}{|l|l|l|l|l|l|l|l|l|r|} \hline \hline
              &\vline \, 1A  &\vline \, 2A  &\vline \, 3A  &\vline \, 3B  & \vline \, 4A
&\vline \, 5A & \vline \, 5B
\\ \hline \hline

$\chi_1$ \, &\vline \, 1  &\vline \, 1 &\vline \, 1 &\vline \, 1  & \vline \, 1
&\vline \, 1 &\vline \, 1
\\ \hline \hline

$\chi_2$ \, &\vline \, 5  &\vline \, 1 &\vline \, 2   &\vline \, -1  & \vline \, -1
&\vline \, 0  &\vline \, 0
\\ \hline \hline

$\chi_3$  &\vline \, 5    &\vline \, 1  & \vline \, -1  &\vline \, 2 &\vline \, -1
&\vline \,  0  &\vline \, 0
\\ \hline \hline

$\chi_4$  &\vline \, 8    &\vline \, 0  &\vline \, -1   &\vline \, -1
&\vline \, 0 &\vline \,  $({-1-\sqrt{5}})/2$  &\vline \, $({-1+\sqrt{5}})/{2}$
\\ \hline \hline

$\chi_5$  &\vline \, 8    &\vline \, 0  &\vline \, -1 &\vline \, -1
& \vline \, 0 &\vline \,  $({-1+\sqrt{5}})/{2}$  &\vline \, $({-1-\sqrt{5}})/{2}$
\\ \hline \hline

$\chi_6$  &\vline \, 9    &\vline \ 1   &\vline \, 0  &\vline \, 0
& \vline \, 1 &\vline \, -1  &\vline \, -1
\\ \hline \hline

$\chi_7$  &\vline \, 10   &\vline \ -2  &\vline \, 1  &\vline \, 1  &\vline \, 0
&\vline \,  0  &\vline \, 0
\\ \hline \hline

\end{tabular}

\par \vskip 1pc

\begin{claim} \label{claim:irred} As $A_{6}$-modules, one has
$$S(X) \otimes \mathbf C = \chi_{1} \oplus \chi_{2} \oplus \chi_{3} \oplus \chi_{6}\,\, .$$
\end{claim}

\begin{proof} Since $S(X)^{A_{6}} = \mathbf Z H$,  the irreducible decomposition must be of the
following form:

$$S(X) \otimes \mathbf C = \chi_{1} \oplus a_{2}\chi_{2} \oplus a_{3}\chi_{3} \oplus
a_{4}\chi_{4} \oplus a_{5}\chi_{5} \oplus a_{6}\chi_{6} \oplus
a_{7}\chi_{7}\,\, ,$$ where $a_{i}$ are non-negative integers. Let
us determine $a_{i}$'s. As in (2), using the topological Lefschetz
fixed point formula and the fact that $\text{rank}\, T(X) = 2$, we
have
$$\chi_{\text{top}}(X^{g}) = 4 + \text{tr}(g^{*} \vert S(X))$$
for $g \in A_{6}$. Running $g$ through the $7$-conjugacy classes and
calculating both sides based on Nikulin's table and the character table
above, we obtain the following system of equations:

$$20 = 1 + 5(a_{2} + a_{3}) + 8(a_{4} + a_{5}) + 9a_{6} + 10a_{7}\,\, ,$$
$$4 = 1 +  (a_{2} + a_{3}) + a_{6} - 2a_{7}\,\, ,$$
$$2 = 1 +  (2a_{2} - a_{3}) - (a_{4} + a_{5}) + a_{7}\,\, ,\,\,
2 = 1 +  (-a_{2} +2a_{3}) - (a_{4} + a_{5}) + a_{7}\,\,  ,$$
$$0 = 1 -  (a_{2} + a_{3}) + a_{6}\,\, ,$$
$$0 = 1 -  (\frac{1 + \sqrt{5}}{2}a_{4} + \frac{1 - \sqrt{5}}{2} a_{5}) - a_{6}\,\, ,\,\,
0 = 1 - (\frac{1 - \sqrt{5}}{2}a_{4} + \frac{1 + \sqrt{5}}{2} a_{5})  - a_{6}\,\,  .$$

Now, we get the result by solving this system of Diophantine equations.
\end{proof}

From now, assuming to the contrary
that $I = 3$ and $G_{N} = A_{6}$, we shall derive a contradiction.

\begin{claim} \label{claim:direct} $G = G_{N} \times \mu_{3} (= A_{6} \times \mu_{3})$.
\end{claim}

\begin{proof}  We shall use the following general
proposition by [IOZ] (see also [Og, Proposition (5.1)]):

\begin{proposition} \label{proposition:split}
Assume that $I = 3$. Let $g$ be an element of $G$ such that
$\alpha(g) = \zeta_{3}$, i.e. $g^{*}\omega_{X} =
\zeta_{3}\omega_{X}$. Set $\text{ord}(g) = 3k$. Then $(k, 3) = 1$.
In particular, $G = G_{N} : \mu_{3}$, a semi-direct product.
\end{proposition}

From this proposition and our assumption, we have $G = A_{6} :
\mu_{3}$, a semi-direct product. Let $h$ be an element of $G$ such
that $\alpha(h) = \zeta_{3}$ and $\text{ord}(h) = 3$. Since
$\text{Out}(A_{6}) = C_2^{\oplus {2}}$, it follows that there is
an element $a \in A_{6}$ such that $h^{-1}xh = a^{-1}xa$ for all
$x \in A_{6}$. Then $ha^{-1} \in Z(G)$, $\alpha(ha^{-1}) =
\zeta_{3}$, and $\text{ord}(ha^{-1}) = 3l$ with $(l, 3) = 1$. Here
$Z(G)$ is the center of $G$. So, replacing $ha^{-1}$ by
$(ha^{-1})^{\pm l}$, one obtains an element $g$ such that
$\alpha(g) = \zeta_{3}$, $\text{ord}(g) = 3$ and $gx = xg$ for all
$x \in G$. This implies the result.
\end{proof}

Let $g$ be a generator of $\mu_{3}$ in $G = A_{6} \times \mu_{3}$.
Then $g^{*}$ makes the irreducible decomposition in Claim (4.2)
stable, i.e. $g^{*} (\chi_{i}) = \chi_{i}$. By the Schur lemma, $g
\vert \chi_{i}$ is a scalar multiplication. Moreover, by $g^{*}H =
H$, one has $g^{*} \vert \chi_{1} = 1$. Set $g^{*} \vert \chi_{2}
= \zeta_{3}^{a}$, $g^{*} \vert \chi_{3} = \zeta_{3}^{b}$ and
$g^{*} \vert \chi_{6} = \zeta_{3}^{c}$, where $a, b, c \in \{0, 1,
2\}$. Note that the action of $g^{*}$ on $S(X)\otimes \mathbf C$
is defined over $\mathbf Z$. Thus, the multiplicities of
$\zeta_{3}$ and $\zeta_{3}^{2} (= \zeta_{3}^{-1})$ must be equal.
So by $\dim \chi_{2} = \dim \chi_{3} = 5$ and $\dim \chi_{6} = 9$,
one has $c = 0$ and $a + b \equiv 0\, \text{mod}\, 3$, i.e. $(a,
b, c) = (0, 0, 0)$, $(1, 2, 0)$ or $(2, 1, 0)$.

Let us first consider the case $(a, b, c) = (0,0,0)$. In this
case, we have $g^{*} \vert S(X) = id$. Moreover, we have
$\text{rank}\, T(X) = 2 = \varphi(3)$. Here $\varphi$ is the Euler
function. Thus, by [OZ1, main Theorem], the intersection matrix of
$T(X)$ is
$$\begin{pmatrix}
2&1\cr
1&2
\end{pmatrix}\,\, .$$
Therefore by [SI], one has
$$X \simeq\,\,   \text{the\, minimal\, resolution\, of}\,
(E_{\zeta_{3}} \times  E_{\zeta_{3}})/\langle \text{diag}\, (\zeta_{3}, \zeta_{3}^{2}) \rangle\,\, .$$
So by the main result of Vinberg [Vi], we have
$$\text{Aut}(X) \simeq C_{3} \times (C_{2}^{* 12} :((S_{3} \times S_{3}):C_{2}))\,\,  ,$$
where $C_{2}^{* 12}$ denotes the free product of $12$ $C_{2}$'s. Since
$A_{6} (< \text{Aut}(X))$ is simple
and $C_{3} < \text{Aut}\, X$ is normal, we have $A_{6} \cap C_{3} = \{1\}$. Hence $A_{6}$ becomes a subgroup
of the quotient group $(C_{2}^{* 12}:((S_{3} \times S_{3}):C_{2}))$. Since $C_{2}^{*12}$ has
no elements of finite order, except involutions and identity, we have again
$A_{6} \cap C_{2}^{*12} = \{1\}$ for the same reason, and $A_{6}$ becomes a subgroup
of the quotient group $(S_{3} \times S_{3}):C_{2}$. However, $\vert (S_{3} \times S_{3}):C_{2}
\vert = 72$, while $\vert A_{6} \vert = 360$, a contradiction.

Let us consider the case $(a, b, c) = (1, 2, 0)$ (resp. $(a, b, c)
= (2, 1, 0)$). Let $\tau \in A_{6}$ be an element of the conjugacy
class $3A$ (resp. $3B$).

By Claim (4.2) and by the case assumption together with the
character table, we have
$$\text{tr}\, (\tau g)^{*} \vert S(X) = 1 +  2\zeta_{3} - \zeta_{3}^{2}\,,$$
which is a contradiction, as the left hand side is an integer
while the right hand side is not. This completes the proof of (3) and
also Proposition (4.1).
\end{proof}

The next proposition is the main result of this section.

\begin{proposition} \label{proposition:uniquesurf}
Assume that $X$ is a $K3$ surface admitting a faithful group action
of $A_{6}.\mu_{4}$. Set $G = A_{6}.\mu_{4}$. Then
\begin{list}{}{
\setlength{\leftmargin}{10pt}
\setlength{\labelwidth}{6pt}
}
\item[(1)] $X \simeq F$, where $F$ is the surface constructed in Theorem (3.1).
\item[(2)] $S(X)^{G} = S(X)^{G_{N}}
= \mathbf Z H$ where $H$ is an ample primitive class with
$(H^{2}) = 20$.
\end{list}
\end{proposition}

Let us first show the following:
\begin{proposition}\label{proposition:rank}
$A_{L^{G_{N}}} := (L^{G_{N}})^{*}/L^{G_{N}} \simeq \mathbf Z/3 \oplus \mathbf Z/60$.
\end{proposition}

\begin{proof}
We denote by $N(Rt)$ the (non-Leech) Niemeier lattice whose root
lattice is isomorphic to $Rt$. There are $23$ such lattices up to
isomorphisms [CS, Chapter 18].  By [Ko2, Lemmas 5 and 6], there
are a (non-Leech) Niemeier lattice $N(Rt)$, a primitive embedding
$A_{1} \oplus L_{G_{N}} \subset N(Rt)$ and a faithful action of
$G_{N}$ on $N(Rt)$ such that $G_{N}$ acts on $A_{1}$ trivially.
Moreover, one can choose this action so that $L_{G_{N}} =
N(Rt)_{G_{N}}$, and a Weyl chamber (one of whose codimension one
faces corresponds to $A_{1}$) stable. Here $N(Rt)_{G_{N}}$ is the
orthogonal complement of the sublattice $N(Rt)^{G_{N}}$ in
$N(Rt)$. So $G_N$ can be regarded as a subgroup of the symmetry
group of the Coxeter-Dynkin diagram of $Rt$, i.e. one has
$$G_{N}
< S(N(Rt)) := O(N(Rt))/W(N(Rt)),$$
 where $W(N(Rt))$ is the Weyl
group. Furthermore, $G_N$ fixes a vertex of the diagram.

Note that $\text{rank}\, N(Rt)^{G_{N}} = \text{rank}\, L^{G_{N}} +
2 = 5$ and $ A_{L^{G_{N}}} \simeq A_{L_{G_{N}}}(-1) =
A_{N(Rt)_{G_{N}}}(-1) \simeq A_{N(Rt)^{G_{N}}}$. Now Proposition
(4.6) follows from Lemma (4.7) below.
\end{proof}

\begin{lemma} \label{lemma:Niem} \begin{list}{}{
\setlength{\leftmargin}{10pt}
\setlength{\labelwidth}{6pt}
}
\item[(1)] The root lattice $Rt$ is either $A_{1}^{\oplus 24}$ or $A_{2}^{\oplus 12}$.
\item[(2)] Assume that the first case in (1) occurs. Then the orbit decomposition of the action of $G_{N}$ on the $24$ simple roots
is either $[1, 1, 1, 6, 15]$ or $[1,1,6,6,10]$. The intersection matrix
of $N(A_{1}^{\oplus 24})^{G_{N}}$ (under some integral basis) is then
either
$$\begin{pmatrix}
-2&0&-1&0&0\cr
0&-2&-1&0&0\cr
-1&-1&-4&0&0\cr
0&0&0&-2&-1\cr
0&0&0&-1&-8\cr
\end{pmatrix}\,\, \text{or}\,\, \begin{pmatrix}
-2&0&-1&-1&-1\cr
0&-2&-1&-1&-1\cr
-1&-1&-4&-1&-1\cr
-1&-1&-1&-4&-1\cr
-1&-1&-1&-1&-6\cr
\end{pmatrix}\,\, .$$
In particular, $A_{N(Rt)^{G_{N}}} \simeq \mathbf Z/3 \oplus \mathbf Z/60$.
\item[(3)] Assume that the second case in (1) occurs. Then the orbit decomposition of the action of $G_{N}$ on the $24$ simple roots
is $[1, 1, 1, 1, 20]$. The intersection matrix
of $N(A_{2}^{\oplus 12})^{G_{N}}$ (under some integral basis) is
$$\begin{pmatrix}
-2&1&0&0&0\cr
1&-2&0&0&0\cr
0&0&-2&1&0\cr
0&0&1&-2&0\cr
0&0&0&0&-20\cr
\end{pmatrix}\,\, .$$
In particular, $A_{N(Rt)^{G_{N}}} \simeq \mathbf Z/3 \oplus \mathbf Z/60$
as well.

\end{list}
\end{lemma}

\begin{proof} By the description of $S(N(Rt))$ ([CS, Chapter 18], [Ko2, the proof of
Main Theorem]), and by the fact that $A_{6} (< S(N(Rt)))$ fixes a
vertex of the Coxeter-Dynkin diagram of $Rt$, we see that $Rt$ is
either $A_{1}^{\oplus 24}$ or $A_{2}^{\oplus 12}$. Let $$\mathcal
N_{2}:= \{r_{1}\, ,\, r_{2}\, ,\, \cdots \, , r_{24}\}$$ be the
set of the simple roots corresponding to the codimension one faces
of the stable Weyl chamber. We may assume that $r_{1}$ is fixed by
$A_{6}$. We shall consider two types of $N(Rt)$ one by one.

First consider the case $N := N(A_{1}^{\oplus 24})$.

\begin{claim} \label{claim:orbits} The orbit decomposition type of
$A_{6}$ on $\mathcal N_{2}$ is either (i) $[1, 1, 1, 6, 15]$ or
(ii) $[1, 1, 6, 6, 10]$.  (Note in particular that $A_{6} <
M_{23}$ under the action on $\mathcal N_{2}$.)
\end{claim}

\begin{proof} Since $\text{rank}\, N^{A_{6}} =  5$,
$\mathcal N_{2}$ is divided into exactly 5 orbits, and one of the
orbits is a one-point set, say $[1, a, b, c, d]$ with $1 + a + b +
c + d = 24$. Since $A_{6}$ is a simple group, the action on each
orbit is faithful unless it is a one-point orbit. Thus, $a \vert
360$, and $a \geq 6$ unless $a = 1$. The same holds for $b$, $c$,
$d$. Now assuming $a \leq b \leq c \leq d$, we see that the
orbit decompositions is one of the two in (4.8) and three below:
(iii) $[1, 1, 6, 8, 8]$, (iv) $[1, 1, 1, 9, 12]$, (v) $[1, 1, 1,
1, 20]$. In Case (iii) (resp. (iv)), the transitivity of
$A_6$-action on the length-8 (resp. length-12) orbit implies that
$A_6$ has a (stabilizer) subgroup of order $|A_6|/8 = 45$ (resp.
$|A_6|/12 = 30$), which is impossible by utilizing the list of
maximal subgroups of $A_6$ in the Atlas. If Case (v) occurs, then
$A_6 < M_{24}$ fixes 4 letters and hence $A_6 < M_{20}$ so that
$360 = |A_6|$ divides $|M_{20}| = 960$, absurd. This proves the
claim.
\end{proof}

By this claim, after re-numbering the elements of $\mathcal N_{2}$, we have
$$Rt^{A_{6}} = \mathbf Z \langle s_{1}, s_{2}, s_{3}, s_{4}, s_{5} \rangle\,\, ,$$
where in Case (i) of (4.8)
$$s_{1} = r_{1}\, ,\, s_{2} = r_{2}\, ,\, s_{3} = r_{3}\, ,\,
s_{4} = r_{4} + \cdots + r_{9}\, ,\, s_{5} = r_{10} + \cdots +
r_{24}\, ,$$ and in Case (ii) of (4.8)
$$s_{1} = r_{1}\, ,\, s_{2} = r_{2}\, ,\, s_{3} = r_{3} + \cdots + r_{8}\, ,\,
s_{4} = r_{9} + \cdots + r_{14}\, ,\,
s_{5} = r_{15} + \cdots +  r_{24}\, .$$

\begin{claim} \label{claim:inv} According to the cases (i) and (ii), one has:
$$N^{A_{6}} = \mathbf Z \langle s_{1}, s_{2}, s_{3}, \frac{s_{1} + s_{2} +
s_{4}}{2},
\frac{s_{3} + s_{5}}{2} \rangle \,\,  ,$$
$$N^{A_{6}} = \mathbf Z \langle s_{1}, s_{2},
\frac{s_{1} + s_{2} + s_{3}}{2},  \frac{s_{1} + s_{2} + s_{4}}{2},
\frac{s_{1} + s_{2} + s_{5}}{2} \rangle \,\,  .$$
\end{claim}

\begin{proof} Recall that the subset $N/Rt(\simeq \mathbf F_{2}^{\oplus 12})$ of $Rt^{*}/Rt \simeq \mathbf F_{2}^{\oplus 24}$
is the so-called binary Golay code. The element of length $8$
(resp. $12$) is  called an octad (resp. a dodecad). We often
identify $\sum_{k \in K}r_{k}/2$ with the subset $\{r_{k} \vert k
\in K\}$ of $\mathcal N_{2}$. The set of octads forms the Steiner
system $St(5, 8, 24)$ of $\mathcal N_{2}$. Note that our numbering
of elements is different from the one in Section 2. However, since
our proof does not involve calculations based on Todd's list, we
continue to keep our numbering of elements from $1$ to $24$ (not
from $\infty$, $0$, to $22$).

Recall that $Rt^{A_{6}} \subset N^{A_{6}} \subset
(Rt^{*})^{A_{6}}$, and that the lattice $N^{A_{6}}$ is generated
by $Rt^{A_{6}}$, $\sum_{k=1}^{24}r_{k}/2$ and $\sum_{k \in K}
r_{k}/2$, where $K$ runs through all $A_{6}$-invariant octads,
dodecads, or the complements of octads. In what follows, we
consider the second case, i.e. the case where the orbit
decomposition is
$$\mathcal N_{2} = \mathcal O_{1} \cup \mathcal O_{2} \cup \mathcal O_{3} \cup \mathcal O_{4} \cup \mathcal O_{5}\,\, ,$$
where
$$\vert \mathcal O_{1} \vert = \vert \mathcal O_{2} \vert = 1\,\, ,\,\,
 \vert \mathcal O_{3} \vert = \vert \mathcal O_{4} \vert = 6\,\, ,\,\,
\vert \mathcal O_{5} \vert = 10\,\, .$$ The first case is easier,
so we omit its proof.

By the shape of the orbit decomposition, the candidates of
$A_{6}$-invariant octads and dodecads are
$$(s_{1} + s_{2} + s_{3})/2\,  ,\,  (s_{1} + s_{2} + s_{4})/2\,  ,\,
(s_{1} + s_{2} + s_{5})/2\,  ,\,   (s_{3} + s_{4})/2\,  .$$ Let us
show that the first sum, or equivalently, the set $\mathcal O_{1}
\cup \mathcal O_{2} \cup \mathcal O_{3}$, is indeed an octad.
Choose an order $5$ element $g \in A_{6} (< M_{23})$. Then, by
[EDM, Appendix B, Table 5.I],  the cycle decomposition type of $g$
on $\mathcal N_{2}$ is $(1^{4})(5^{4})$. Thus, the cycle type of
the action of $g$ on $\mathcal O_{3}$ is $(1)(5)$. So, after
re-numbering elements in $\mathcal O_{3}$, we may assume that
$g(r_{i}) = r_{i+1}$ ($3 \leq i \leq 6$), $g(r_{7}) = r_{3}$ and
$g(r_{8}) = r_{8}$. By the Steiner property, there is an octad $A$
containing the $5$-element set $\{r_{3}, r_{4}, r_{5}, r_{6},
r_{7}\}$. Since
$$g(\{r_{3}, r_{4}, r_{5}, r_{6}, r_{7}\})
= \{r_{3}, r_{4}, r_{5}, r_{6}, r_{7}\}$$ we have by the Steiner
property that $A = g(A)$. Let $s$, $t$, $u$ be the remaining three
elements of $A$. Since $g$ acts on the fifth orbit $\mathcal
O_{5}$ as $(5^{2})$, none of them is in the fifth orbit. If two or
three of $s$, $t$, $u$ are in the fourth orbit $\mathcal O_{4}$,
using the fact that the cycle type of $g$ on the fourth orbit is of
(1)(5), we have $g^{m}(A) \not= A$ for a suitable $m$, a
contradiction. Assume that exactly one of $s$, $t$, $u$, say $s$,
is in the fourth orbit $\mathcal O_{4}$.  Note that $\mathcal
O_{4}$ and $\mathcal O_{3}$ are both order $6$ set. Then the cycle
type of an involution $\tau \in A_{6}$, is necessarily of type
$(1^{2})(2^{2})$ both on $\mathcal O_{4}$ and on $\mathcal O_{3}$.
So, there is an involution $\tau \in A_{6}$ such that $\tau(s)
\not\in A$, but at least one of $\tau(t)$, $\tau(u)$, and at least
four of $\tau(r_{3})\, , \cdots \,  , \tau(r_{7})$ are in $A$.
Thus, $\vert \tau(A) \cap A \vert \geq 5$, whence $\tau(A) = A$,
by the Steiner property, a contradiction to $s \in A$ but $\tau(s)
\not\in A$. So, none of $s$, $t$, $u$ is in $\mathcal O_{4} \cup
\mathcal O_{5}$. This means that $A$ is the union of the first
three orbits and is then $A_{6}$-invariant. In the exactly same
manner,  the union of the first two and the fourth orbits is also
an $A_{6}$-invariant octad. Then, by taking a symmetric difference
and complement, we find that the other two candidates are actually
$A_{6}$-invariant dodecads. This implies the result.
\end{proof}

The matrices in Lemma (4.7)(2) are nothing but the intersection matrices
with respect to these basis in (4.9). Calculating elementary divisors, we also get the last statement of Lemma (4.7)(2).

\par \vskip 1pc
Next consider the case $N := N(A_{2}^{\oplus 12})$.

\begin{claim} \label{claim:orbit} The orbit decomposition type of
$A_{6}$ on $\mathcal N_{2}$ is $[1, 1, 1, 1, 20]$.
\end{claim}

\begin{proof} As before,
$\mathcal N_{2}$ is divided into exactly 5 orbits, and one of the
orbits is a one point set. Since the graph of $\mathcal N_{2}$
consists of $12$-connected components,

\par \vskip 3pc


\setlength{\unitlength}{3947sp}%
\begingroup\makeatletter\ifx\SetFigFont\undefined%
\gdef\SetFigFont#1#2#3#4#5{%
  \reset@font\fontsize{#1}{#2pt}%
  \fontfamily{#3}\fontseries{#4}\fontshape{#5}%
  \selectfont}%
\fi\endgroup%
\begin{picture}(5333,3553)(301,-3044)
\thinlines
{ \put(5476, 89){\vector( 0,-1){0}}
\put(3264, 89){\oval(4424,300)[tr]}
\put(3264, 89){\oval(4426,300)[tl]}
\put(1051, 89){\vector( 0,-1){0}}
}%
{ \put(1576,-1861){\vector( 0, 1){0}}
\put(3300,-1861){\oval(3448,306)[bl]}
\put(3300,-1786){\oval(4202,456)[br]}
\put(5401,-1786){\vector( 0, 1){0}}
}%
{ \put(601,-361){\circle{150}}
}%
{ \put(601,-1261){\circle{150}}
}%
{ \put(1051,-361){\circle{150}}
}%
{ \put(1051,-1261){\circle{150}}
}%
{ \put(1501,-361){\circle{150}}
}%
{ \put(1501,-1261){\circle{150}}
}%
{ \put(1951,-361){\circle{150}}
}%
{ \put(1951,-1261){\circle{150}}
}%
{ \put(2401,-361){\circle{150}}
}%
{ \put(2401,-1261){\circle{150}}
}%
{ \put(2851,-361){\circle{150}}
}%
{ \put(2851,-1261){\circle{150}}
}%
{ \put(3301,-361){\circle{150}}
}%
{ \put(3301,-1261){\circle{150}}
}%
{ \put(3751,-361){\circle{150}}
}%
{ \put(3751,-1261){\circle{150}}
}%
{ \put(4201,-361){\circle{150}}
}%
{ \put(4201,-1261){\circle{150}}
}%
{ \put(4651,-361){\circle{150}}
}%
{ \put(4651,-1261){\circle{150}}
}%
{ \put(5101,-361){\circle{150}}
}%
{ \put(5101,-1261){\circle{150}}
}%
{ \put(5551,-361){\circle{150}}
}%
{ \put(5551,-1261){\circle{150}}
}%
{ \put(601,-436){\line( 0,-1){750}}
}%
{ \put(1051,-436){\line( 0,-1){750}}
}%
{ \put(1501,-436){\line( 0,-1){750}}
}%
{ \put(1951,-436){\line( 0,-1){750}}
}%
{ \put(2401,-436){\line( 0,-1){750}}
}%
{ \put(2851,-436){\line( 0,-1){750}}
}%
{ \put(3301,-436){\line( 0,-1){750}}
}%
{ \put(3751,-436){\line( 0,-1){750}}
}%
{ \put(4201,-436){\line( 0,-1){750}}
}%
{ \put(4651,-436){\line( 0,-1){750}}
}%
{ \put(5101,-436){\line( 0,-1){750}}
}%
{ \put(5551,-436){\line( 0,-1){750}}
}%
\put(301,-136){\makebox(0,0)[lb]{\smash{\SetFigFont{17}{20.4}{\rmdefault}{\bfdefault}{\updefault}{ r}%
}}}
\put(301,-1561){\makebox(0,0)[lb]{\smash{\SetFigFont{17}{20.4}{\rmdefault}{\bfdefault}{\updefault}{ r}%
}}}
\put(451,-211){\makebox(0,0)[lb]{\smash{\SetFigFont{14}{16.8}{\rmdefault}{\mddefault}{\updefault}{ 1 }%
}}}
\put(751,-136){\makebox(0,0)[lb]{\smash{\SetFigFont{17}{20.4}{\rmdefault}{\bfdefault}{\updefault}{ r}%
}}}
\put(901,-211){\makebox(0,0)[lb]{\smash{\SetFigFont{14}{16.8}{\rmdefault}{\mddefault}{\updefault}{ 3}%
}}}
\put(1201,-136){\makebox(0,0)[lb]{\smash{\SetFigFont{17}{20.4}{\rmdefault}{\bfdefault}{\updefault}{ r}%
}}}
\put(1351,-211){\makebox(0,0)[lb]{\smash{\SetFigFont{14}{16.8}{\rmdefault}{\mddefault}{\updefault}{ 5}%
}}}
\put(2176,-136){\makebox(0,0)[lb]{\smash{\SetFigFont{17}{20.4}{\rmdefault}{\bfdefault}{\updefault}{ r}%
}}}
\put(2326,-211){\makebox(0,0)[lb]{\smash{\SetFigFont{14}{16.8}{\rmdefault}{\mddefault}{\updefault}{ 9}%
}}}
\put(2551,-136){\makebox(0,0)[lb]{\smash{\SetFigFont{17}{20.4}{\rmdefault}{\bfdefault}{\updefault}{ r}%
}}}
\put(2701,-211){\makebox(0,0)[lb]{\smash{\SetFigFont{14}{16.8}{\rmdefault}{\mddefault}{\updefault}{ 11}%
}}}
\put(3001,-136){\makebox(0,0)[lb]{\smash{\SetFigFont{17}{20.4}{\rmdefault}{\bfdefault}{\updefault}{ r}%
}}}
\put(3151,-211){\makebox(0,0)[lb]{\smash{\SetFigFont{14}{16.8}{\rmdefault}{\mddefault}{\updefault}{ 13}%
}}}
\put(3601,-211){\makebox(0,0)[lb]{\smash{\SetFigFont{14}{16.8}{\rmdefault}{\mddefault}{\updefault}{ 15}%
}}}
\put(3451,-136){\makebox(0,0)[lb]{\smash{\SetFigFont{17}{20.4}{\rmdefault}{\bfdefault}{\updefault}{ r}%
}}}
\put(4051,-211){\makebox(0,0)[lb]{\smash{\SetFigFont{14}{16.8}{\rmdefault}{\mddefault}{\updefault}{ 17}%
}}}
\put(3901,-136){\makebox(0,0)[lb]{\smash{\SetFigFont{17}{20.4}{\rmdefault}{\bfdefault}{\updefault}{ r}%
}}}
\put(4351,-136){\makebox(0,0)[lb]{\smash{\SetFigFont{17}{20.4}{\rmdefault}{\bfdefault}{\updefault}{ r}%
}}}
\put(4501,-211){\makebox(0,0)[lb]{\smash{\SetFigFont{14}{16.8}{\rmdefault}{\mddefault}{\updefault}{ 19}%
}}}
\put(4801,-136){\makebox(0,0)[lb]{\smash{\SetFigFont{17}{20.4}{\rmdefault}{\bfdefault}{\updefault}{ r}%
}}}
\put(4951,-211){\makebox(0,0)[lb]{\smash{\SetFigFont{14}{16.8}{\rmdefault}{\mddefault}{\updefault}{ 21}%
}}}
\put(5251,-136){\makebox(0,0)[lb]{\smash{\SetFigFont{17}{20.4}{\rmdefault}{\bfdefault}{\updefault}{ r}%
}}}
\put(5401,-211){\makebox(0,0)[lb]{\smash{\SetFigFont{14}{16.8}{\rmdefault}{\mddefault}{\updefault}{ 23}%
}}}
\put(751,-1561){\makebox(0,0)[lb]{\smash{\SetFigFont{17}{20.4}{\rmdefault}{\bfdefault}{\updefault}{ r}%
}}}
\put(901,-1636){\makebox(0,0)[lb]{\smash{\SetFigFont{14}{16.8}{\rmdefault}{\mddefault}{\updefault}{ 4}%
}}}
\put(451,-1636){\makebox(0,0)[lb]{\smash{\SetFigFont{14}{16.8}{\rmdefault}{\mddefault}{\updefault}{ 2}%
}}}
\put(3226,-1636){\makebox(0,0)[lb]{\smash{\SetFigFont{14}{16.8}{\rmdefault}{\mddefault}{\updefault}{ 14}%
}}}
\put(3526,-1561){\makebox(0,0)[lb]{\smash{\SetFigFont{17}{20.4}{\rmdefault}{\bfdefault}{\updefault}{ r}%
}}}
\put(3676,-1636){\makebox(0,0)[lb]{\smash{\SetFigFont{14}{16.8}{\rmdefault}{\mddefault}{\updefault}{ 16}%
}}}
\put(3151,314){\makebox(0,0)[lb]{\smash{\SetFigFont{17}{20.4}{\rmdefault}{\bfdefault}{\updefault}{ 11}%
}}}
\put(3376,-2386){\makebox(0,0)[lb]{\smash{\SetFigFont{17}{20.4}{\rmdefault}{\bfdefault}{\updefault}{ 10}%
}}}
\put(3076,-1561){\makebox(0,0)[lb]{\smash{\SetFigFont{17}{20.4}{\rmdefault}{\bfdefault}{\updefault}{ r}%
}}}
\put(2101,-1561){\makebox(0,0)[lb]{\smash{\SetFigFont{17}{20.4}{\rmdefault}{\bfdefault}{\updefault}{ r}%
}}}
\put(2251,-1636){\makebox(0,0)[lb]{\smash{\SetFigFont{14}{16.8}{\rmdefault}{\mddefault}{\updefault}{ 10}%
}}}
\put(3976,-1561){\makebox(0,0)[lb]{\smash{\SetFigFont{17}{20.4}{\rmdefault}{\bfdefault}{\updefault}{ r}%
}}}
\put(4126,-1636){\makebox(0,0)[lb]{\smash{\SetFigFont{14}{16.8}{\rmdefault}{\mddefault}{\updefault}{ 18}%
}}}
\put(4426,-1561){\makebox(0,0)[lb]{\smash{\SetFigFont{17}{20.4}{\rmdefault}{\bfdefault}{\updefault}{ r}%
}}}
\put(4576,-1636){\makebox(0,0)[lb]{\smash{\SetFigFont{14}{16.8}{\rmdefault}{\mddefault}{\updefault}{ 20}%
}}}
\put(4876,-1561){\makebox(0,0)[lb]{\smash{\SetFigFont{17}{20.4}{\rmdefault}{\bfdefault}{\updefault}{ r}%
}}}
\put(5026,-1636){\makebox(0,0)[lb]{\smash{\SetFigFont{14}{16.8}{\rmdefault}{\mddefault}{\updefault}{ 22}%
}}}
\put(5326,-1561){\makebox(0,0)[lb]{\smash{\SetFigFont{17}{20.4}{\rmdefault}{\bfdefault}{\updefault}{ r}%
}}}
\put(5476,-1636){\makebox(0,0)[lb]{\smash{\SetFigFont{14}{16.8}{\rmdefault}{\mddefault}{\updefault}{ 24}%
}}}
\put(2626,-1561){\makebox(0,0)[lb]{\smash{\SetFigFont{17}{20.4}{\rmdefault}{\bfdefault}{\updefault}{ r}%
}}}
\put(2776,-1636){\makebox(0,0)[lb]{\smash{\SetFigFont{14}{16.8}{\rmdefault}{\mddefault}{\updefault}{ 12}%
}}}
\put(1651,-136){\makebox(0,0)[lb]{\smash{\SetFigFont{17}{20.4}{\rmdefault}{\bfdefault}{\updefault}{ r}%
}}}
\put(1876,-211){\makebox(0,0)[lb]{\smash{\SetFigFont{14}{16.8}{\rmdefault}{\mddefault}{\updefault}{ 7}%
}}}
\put(1201,-1561){\makebox(0,0)[lb]{\smash{\SetFigFont{17}{20.4}{\rmdefault}{\bfdefault}{\updefault}{ r}%
}}}
\put(1351,-1636){\makebox(0,0)[lb]{\smash{\SetFigFont{14}{16.8}{\rmdefault}{\mddefault}{\updefault}{ 6}%
}}}
\put(1651,-1561){\makebox(0,0)[lb]{\smash{\SetFigFont{17}{20.4}{\rmdefault}{\bfdefault}{\updefault}{ r}%
}}}
\put(1801,-1636){\makebox(0,0)[lb]{\smash{\SetFigFont{14}{16.8}{\rmdefault}{\mddefault}{\updefault}{ 8}%
}}}
\put(2701,-2986){\makebox(0,0)[lb]{\smash{\SetFigFont{14}{16.8}{\rmdefault}{\mddefault}{\updefault}{ Figure 6}%
}}}
\end{picture}

\vskip 1pc \noindent
the orbit decomposition is $[1, 1, a, b, c]$ with $a +
b + c = 22$ ($a \leq b \leq c$) and $A_{6}$ acts on the set of $11$-connected
components. Since $A_{6}$ is simple,  if $A_{6}$ acts
on $m$ components with $m \leq 5$, then it acts trivially on each
of the $m$ components, whence $m = 1$ and its $2$ simple roots are
both fixed. Thus, $a = b = 1$ and $c = 20$.
\end{proof}

By this claim, after re-numbering the elements of $\mathcal N_{2}$, we have
$$Rt^{A_{6}} = \mathbf Z \langle s_{1}, s_{2}, s_{3}, s_{4}, s_{5} \rangle$$
where $$s_{1} = r_{1}\, ,\, s_{2} = r_{2}\, ,\, s_{3} = r_{3}\,
,\, s_{4} = r_{4}\, ,\, s_{5} = r_{5} + \cdots +  r_{24}\, ,$$ and
$r_{i}$ are labeled so that $\{r_{2k-1}, r_{2k}\}$ forms a
connected component.
\begin{claim} \label{claim:inva} One has:
$$N^{A_{6}} = \mathbf Z \langle s_{1}, s_{2}, s_{3}, s_{4},
s_{5}\rangle \,\,  .$$
\end{claim}
\begin{proof} By [CS, Chapter 18], $N^{A_{6}}$ is generated by $Rt^{A_{6}}$ and $A_{6}$-invariant
elements of the form
$v_{S} := \sum_{k \in S} \pm (r_{2k-1} + 2r_{2k})/3$, where $S$ is an element
of the so called ternary Golay code. Such $v_{S}$ must also satisfy $(v_{S}^{2}) \in \mathbf Z$.
However, by the shape of the orbit decomposition, there is no such $A_{6}$-invariant sum.
\end{proof}

Again, the matrix in Lemma (4.7)(3) is nothing but the intersection matrix
with respect to this basis. Calculating elementary divisors, we also get
the last statement of Lemma (4.7)(3).
This proves Lemma (4.7) and also Proposition (4.6).
\end{proof}

In order to complete the proof of Proposition (4.5), we need two more lemmas.

\begin{lemma} \label{lemma:tpart} There is an integral basis $\langle v_{1},
v_{2} \rangle$ of $T(X)$ such that $G/G_{N} = \langle g \rangle$ acts as
$g(v_{1}) = v_{2}$ and $g(v_{2}) =  -v_{1}$. With respect to this basis
$\langle v_{1}, v_{2}\rangle$,
the intersection matrix of $T(X)$ has the following form:
$$\begin{pmatrix}
2m&0\cr
0&2m
\end{pmatrix}\,\, \text{for\, some}\,\, m \in \mathbf Z\,\, .$$
\end{lemma}
\begin{proof} Since $G/G_{N} \simeq \mu_{4}$, the result follows
from [Ko3, Page 1249] (see also [MO], [OZ2, Lemma 2.8]).
\end{proof}

\begin{lemma} \label{lemma:index} Set $l := [L^{G_{N}} : \mathbf Z H
\oplus T(X)]$, where $G_{N} \simeq A_{6}$ and $H$ is the same as in Proposition (4.1)(1).
Then $l = 1$ or $2$. Moreover, if $l = 2$, then
$$ L^{A_{6}} = \mathbf Z \langle \frac{H + v_{1} + v_{2}}{2}, v_{1}, v_{2}\rangle\,\, .$$
Here $\langle v_{1}, v_{2} \rangle$ is an integral basis of $T(X)$ as in
Lemma (4.12).
\end{lemma}

\begin{proof} The proof is identical to [Ko3, Page 1248], [OZ2, Page 177] or [Og, Lemma 6.10]. \end{proof}

Now we are ready to complete the proof of Proposition (4.5).

\begin{proof} Let $(H^{2}) = 2n$ for some positive
integer $n$. Let $m$ be the positive integer
in Lemma (4.12). By virtue of [SI], it suffices to show that $(n, m) = (10, 3)$.
Let $l$ be the same as in Lemma (4.13).

Assume first that $l = 1$. Then $L^{A_{6}} = \mathbf Z H \oplus T(X)$. However, by Proposition (4.6),
we have then $3 \cdot 60 = 2n \cdot 4m^{2}$, and $nm \not\in \mathbf Z$, a contradiction.

Therefore $ l = 2$ and by Proposition (4.6), we have $4 \cdot 3 \cdot 60 = 2n \cdot 4m^{2}$,
that is, $nm^{2} = 90$. Thus, $(n, m) = (90, 1)$ or $(10, 3)$. However if
$(n, m) = (90, 1)$, then the intersection matrix of $L^{A_{6}}$ with respect to
the integral basis $\langle v_{1}, v_{2}, (v_{1} + v_{2} + H)/2 \rangle$ would be
$$\begin{pmatrix}
2m&0&m\cr
0&2m&m\cr
m&m&m + \frac{n}{2}
\end{pmatrix}\,\, = \begin{pmatrix}
2&0&1\cr
0&2&1\cr
1&1&46
\end{pmatrix}\,\, ,$$
whence $A_{L^{G_{N}}} \simeq \mathbf Z/180$, a contradiction to Proposition (4.6).
Thus $(n, m) = (10, 3)$.
\end{proof}

\section{Uniqueness of the $A_{6}.\mu_{4}$-action}

In this section, we shall show the uniqueness of the $A_{6}.\mu_{4}$-action.
\begin{theorem} \label{theorem:uniquegp} Let $X$
be a $K3$ surface admitting a faithful $A_{6}.\mu_{4}$-action, say,
$$\tau : A_{6}.\mu_{4} \times X \longrightarrow X\, .$$
Then the triplet $(X, A_{6}.\mu_{4}, \tau)$ is isomorphic to the triplet
$(F, \tilde{A}_{6}, \rho_{F})$ constructed in Theorem $(3.1)$.
\end{theorem}

\begin{proof} Our argument here is much inspired by [Ko3].
Put $G := G_{\tau} = A_{6}.\mu_{4}$.  Note
that $G_{N} = A_{6}$. By Proposition (4.5), we may assume that $X
= F$. We put $\tau_{F} = \tau$. We denote by $H_{\tau}$ a
primitive $G_{\tau}$-invariant ample class (now) on $F$. By (4.1),
we have $S(F)^{G_{\tau}} = \mathbf Z H_{\tau}$. As before, we set
$L := H^{2}(F, \mathbf Z)$ and put $T := T(F)$ and $S := S(F)$.

We denote by $H_{\rho}$ a primitive $\tilde{A}_{6}$-invariant
ample class on $F$ (see Theorem 3.1 and Remark 3.4).

Let $\mathcal R := \{c, z, x_{0}, r_{0}, x_{1}, x_{2}\}$ be the
set defined in Section 2. Recall that the pointwise stabilizer
group $\text{Aut}(\mathcal D, \mathcal R)$ is isomorphic to
$A_{6}$, that the set $\mathcal R$ generates the primitive
sublattice $R$, being isomorphic to $A_{2}^{\oplus 2} \oplus
A_{1}^{\oplus 2}$, and that we have the following diagram constructed in the proof of Theorem (3.1):
$$T \subset L \supset S = R^{\perp}_{\Pi} \subset \Pi = \Lambda \oplus U
\supset R\,\, ,$$
under which $H_{\rho} = h$ and the generators $t_{1}/6$, $t_{2}/6$ of
$A_{T(F)}$ correspond to the generators $e_{1}$, $e_{2}$ of $A_{R}$
respectively. See Section 2 for the definition of $e_{1}$ and $e_{2}$.

Consider the group action $\tau_{F, S}^{*} : A_{6}
\longrightarrow \text{GL}(S)$ induced by the geometric action of
$\tau_{F}$ on $F$. The action of $\tau_{F, S}^{*}$ is faithful,
because any action of $A_{6}$ on $F$ is symplectic (see
(4.1)). Since $A_{6}$ is simple and $O(A_{S}, q_{S}) \simeq D_8
\times \mathbf Z/2$ by Proposition (2.6), the natural homomorphism
$\tau_{F, S}^{*}(A_{6}) \longrightarrow O(A_{S}, q_{S})$ is
trivial. Thus, the action $\tau_{F, S}^{*}$ on $S$ can be extended
to the action $\tau_{F, \Pi}$ on $\Pi$ in such a way that
$\tau_{F, \Pi} \vert S = \tau_{F, S}^{*}$ and that $\tau_{F, \Pi}
\vert R = id$.

Let $w_{R}$ be the element of $R^{*}$ defined in Lemma (2.8):
$$w_{R} = -(c + z + r_{0} + x_{0} + \frac{x_{1} + x_{2}}{2})\,\, .$$

\begin{lemma} \label{lemma:leechtoo} Consider the element
$w_{\tau} := H_{\tau}/2 + w_{R}$ of
$\Pi \otimes \mathbf Q$. Then, $w_{\tau}$ is a primitive element of
$\Pi$ and $(w_{\tau}^{2}) = 0$.
\end{lemma}
\begin{proof}
We can choose $t_{1}$ and $t_{2}$ as $v_{1}$ and $v_{2}$ in Lemma (4.12).
Then under the natural isomorphism $A_{T} \simeq A_{R}$ induced by the diagram above, we have $t_{1}/6 \leftrightarrow e_{1}$ and $t_{2}/6 \leftrightarrow e_{2}$, and hence
$$A_{T} \ni (t_{1} + t_{2})/2 \leftrightarrow w_{R} \in A_{R}\,\, .$$
On the other hand, we have also $(H_{\tau} + t_{1} + t_{2})/2 \in L$ by
Lemma (4.13), in which we now know that $l = 2$. Thus
$(t_{1} + t_{2})/2$, $H_{\tau}/2$ and $w_{R}$ give the same element of
$A_{T(F)} = A_{S(F)} = A_{R}$ again under the natural identification induced by the diagram above. Since $\Pi$ is unimodular, this implies $w_{\tau} \in \Pi$.
Since $(w_{\tau}, z) = 1$, it follows that $w_{\tau}$ is also primitive.
We can check $(w_{\tau}^{2}) = 0$ by a direct calculation.
\end{proof}

The two elements $w_{\tau}$ and $z$ of $\Pi$ generate
a sublattice $U_{\tau}$, which is isomorphic to $U$, of $\Pi$. Set
$N := U_{\tau}^{\perp}\subset{\Pi}$. This $N$ is a negative
definite even unimodular lattice of rank $24$
and satisfies $\Pi = N \oplus
U_{\tau}$.
\begin{lemma} \label{lemma:leechtoo} $N$ is isomorphic to the Leech lattice.
\end{lemma}
\begin{proof} Put $K := \tau_{F, \Pi}(A_{6})$. Note that $K \simeq A_{6}$. Since $w_{\tau}$
and $z$ are fixed by $K$, this group $K$ acts on its orthogonal
complement $N$. We have $\Pi^{K} = U_{\tau} \oplus N^{K}$ and
$\Pi_K = N_{K} = S_{K}$ in $\Pi$.  This is because the action of $K$ is
trivial on both $R$ and $U_{\tau}$. Using the unimodularity of
$\Pi$ and Propositions (4.5) and (4.6) and the fact that
$\vert \text{det}\, S_{K} \vert = \vert \text{det}\, L_{K} \vert
= \vert \text{det}\, L^{K} \vert$,
we calculate that
$$\vert \text{det}\, \Pi^{K} \vert = \vert \text{det}\, \Pi_{K} \vert = \vert \text{det}\, S_{K} \vert$$
$$= 180 = \vert \text{det}\, (\mathbf Z H_{\tau} \oplus R) \vert/4\,\, .$$
This shows that $\Pi^{K}$ is generated by $H_{\tau}$, $R$ and
$w_{\tau}$, whence, by $w_{\tau}$ and $R$. On the other hand,
the construction in Theorem (3.1)and Proposition 4.6 tell us that
$\Pi^{\rho_{F, \Pi}(A_{6})}$ is
generated by $w = w_{\rho}$ and $R$. Thus, comparing intersection forms,
one finds that $N^{K}$ is isomorphic to
$\Lambda^{\rho_{F, \Pi}(A_{6})}$. In particular, $N^{K}$ contains
no $-2$ vector.

Assume that $N$ is not isomorphic to the Leech lattice. Then $N$
is isomorphic to some non-Leech Niemeier lattice $N = N(Rt)$.
Then $A_{6} \simeq K < S(N) := O(N)/W(N)$, where $W(N)$ is the Weyl group of
$N = N(Rt)$, and $K$ acts on the set of $24$ simple roots forming $Rt$ [Ko2].
Moreover, $K$ has
exactly $5$-orbits on the set of $24$ simple roots,
because $\text{rank}\, N^{K} = 5$. Let us set
orbit decomposition type by $[a, b, c, d, e]$ with $a \leq b \leq
c \leq d \leq e$. Then $5a \leq a + b + c + d + e = 24$, whence $a
\leq 4$. However, since $A_{6}$ is simple, there is no nontrivial
homomorphism $A_{6} \simeq K \longrightarrow S_{a}$ with $a \leq
4$. Thus $a = 1$, i.e. $K$ fixes at least one root, a contradiction.
Now we are done.
\end{proof}

Therefore, $N \simeq \Lambda$ and there is $h_{1} \in
O(\Pi)$ such that $h_{1}(w_{\tau}) = w$. By Theorem (2.3), there
is $h_{2} \in \text{Aut}(\mathcal D)$ such that
$h_{2}h_{1}(\mathcal R) = \mathcal R$ (and of course
$h_{2}h_{1}(w_{\tau}) = w$). Recall that the natural homomorphism
$\text{Aut}(\mathcal D, \mathcal R) \longrightarrow O(A_{R}, q_{R})$ is
surjective by Proposition (2.6).  So, there is $h_{3} \in
\text{Aut}(\mathcal D, \mathcal R)$ such that $h_{3}h_{2}h_{1}(w_{\tau}) = w$,
$h_{3}h_{2}h_{1}(\mathcal R) = \mathcal R$ and
$h_{3}h_{2}h_{1} \vert \mathcal R = id$.
Set $f = h_{3}h_{2}h_{1}$ and $f_{S} := f \vert
S$. Note that $f_{S} = id$ on $A_{S}$. So, we can extend $f_{S}$
to an isometry $\varphi$ on $L := H^{2}(F, \mathbf Z)$ such that
$\varphi_S := \varphi \vert S = f_S$ and
$\varphi_{T} := \varphi \vert T = id$. Here we put $T := T(F)$. In particular,
$\varphi$ preserves the Hodge decomposition. Recall that $w_{\tau} = H_{\tau}/2 + w_{R}$ by definition (see Lemma (5.2)) and that $(w = ) w_{\rho} = H_{\rho}/2 + w_{R}$ by the construction in Lemma (2.8) (see also Remark 3.4). Then,
$\varphi_{S}(H_{\tau})= H_{\rho}$, and hence $\varphi$ is also
effective. Thus, there is $\psi \in \text{Aut}(F)$ such that
$\psi^{*} = \varphi$ by the global Torelli Theorem.

By the construction of $f$,  we have $f \circ
\tau_{F,\Pi} \circ f^{-1}(A_{6}) \subset
\text{Aut}(\mathcal D, \mathcal R) = A_{6}$, whence $f
\circ \tau_{F,\Pi} \circ f^{-1}(A_{6}) =
\text{Aut}(\mathcal D, \mathcal R)$. On the other hand, by the
construction of $(F, \tilde{A}_{6}, \rho_{F})$ in Theorem (3.1),
we have also $\rho_{F, \Pi}(A_{6}) = \text{Aut}(\mathcal D,
\mathcal R)$. Thus
$$(S, A_{6}, \psi_{S}^{*} \circ \tau_{F, S}^{*} \circ (\psi_{S}^{*})^{-1}) = (S, A_{6}, \rho_{F, S}^{*})\,\, ,$$
up to $\text{Aut}(A_{6})$, that is, there is $\Psi \in
\text{Aut}(A_{6})$,  such that
$$\psi^{*} \circ \tau^*_{F, S}(a)
\circ (\psi^{*})^{-1}(v) = \rho_{F, S}^{*}(\Psi (a))(v)$$
for all $a
\in A_{6}$ and $v \in S$. Since both $G_{N}$ and
$(\tilde{A}_{6})_{N}$ act on $T$ trivially, we have also
$$(T, A_{6}, \psi^{*}_{T} \circ \tau_{F, T}^{*} \circ (\psi_{T}^{*})^{-1}) =
(T, A_{6}, \tau_{F, T}^{*}) = (T, A_{6}, \rho_{F, T}^{*})$$ via the {\it same}
$\Psi \in \text{Aut}(A_{6})$. Hence, we obtain
$$(L, A_{6}, \psi^{*} \circ \tau_{F, L}^{*} \circ (\psi^{*})^{-1}) = (L, A_{6}, \rho_{F, L}^{*})$$
again via $\Psi$. We note that $\psi^{*} \circ \tau_{F, L}^{*}
\circ (\psi^{*})^{-1}$ and $\rho_{F, L}^{*}$ are both defined over
$L$. So, they coincide on $L$ if they coincide on the finite index
sublattice $S \oplus T$ of $L$.  Now, it follows from the
injectivity part of the global Torelli theorem for $K3$ surfaces
that $(F, A_{6}, \psi^{-1} \circ \tau_{F} \circ \psi) \simeq (F,
A_{6}, \rho_{F})$, and we can identify $(F, G_{N}) = (F,
(\tilde{A_{6}})_{N})$. Then $S(F)^{G} = S(F)^{G_{N}} =
S(F)^{(\tilde{A}_{6})_{N}} = S(F)^{\tilde{A_{6}}} = \mathbf Z H$
by Proposition (4.1)(1) applied for $F$. Here we denote by $H =
H_{\rho}$ the $\tilde{A}_{6}$-invariant primitive polarization on
$F$. Thus $G\, , \tilde{A}_{6} < \text{Aut}(F, H)$, i.e. $G$ and
$\tilde{A}_{6}$ are subgroups of the automorphism group of the
polarized $K3$ surface $(F, H)$. Since $\text{Aut}(F, H)$ is a
finite group containing $A_{6}$, it follows from Proposition
(4.1)(2) that $\vert \text{Aut}(F, H) \vert \leq 4\vert A_{6}
\vert$. Since $G, \tilde{A}_{6} \le \text{Aut}(F, H)$ and $\vert G
\vert = \vert \tilde{A}_{6} \vert = 4 \vert A_{6} \vert$, we have
then $G = \text{Aut}(F, H) = \tilde{A}_{6}$. This completes the
proof.
\end{proof}


\end{document}